\documentclass[12pt,reqno,a4paper]{amsart}    

\usepackage[totalwidth=515pt, totalheight=722pt]{geometry}

\usepackage{diagrams}
\usepackage{amsmath, amssymb}
\usepackage{amsthm}

\usepackage{graphicx}                      
\newcommand{\color}[2][{}]{}         

\usepackage{bbm}         

\usepackage{mathrsfs}    
\renewcommand\mathcal\mathscr  

\usepackage{accents}     

\usepackage{slashed}     

\usepackage{stmaryrd}    

\usepackage{upgreek}     

\numberwithin{equation}{section}

\theoremstyle{plain}            
\newtheorem{theorem}{Theorem}[section]
\newtheorem{proposition}[theorem]{Proposition}
\newtheorem{lemma}[theorem]{Lemma}

\theoremstyle{definition}       
\newtheorem{definition}[theorem]{Definition}
\newtheorem{assumption}[theorem]{Assumption}
\newtheorem{example}[theorem]{Example}

\theoremstyle{remark}           
\newtheorem{remark}[theorem]{Remark}

\newcommand{\Sec}[1]{Section~\ref{sec:#1}}

\newcommand{\Eq}[1]{Eq.~\eqref{eq:#1}}

\newcommand{\Fig}[1]{Figure~\ref{fig:#1}}

\newcommand{\Thm}[1]{Theorem~\ref{thm:#1}}
\newcommand{\Thms}[2]{Theorems~\ref{thm:#1} and~\ref{thm:#2}}

\newcommand{\Thmenum}[2]{Theorem~\ref{thm:#1}~(\ref{#2})}

\newcommand{\Exenum}[2]{Example~\ref{ex:#1}~(\ref{#2})}

\newcommand{\Lem}[1]{Lemma~\ref{lem:#1}}

\newcommand{\Prp}[1]{Proposition~\ref{prp:#1}}

\newcommand{\Remenum}[2]{Remark~\ref{rem:#1}~(\ref{#2})}

\newcommand{\card}[1]{\lvert#1\rvert}   

\DeclareMathOperator{\dd}    {d\!} 
\DeclareMathOperator{\dist}   {dist}
\DeclareMathOperator{\dom}    {dom}

\DeclareMathOperator{\id}     {id}  
\DeclareMathOperator{\ind}   {ind}      

\DeclareMathOperator{\ran}    {ran}

\DeclareMathOperator{\vol}    {vol}
\DeclareMathOperator{\tr}     {tr}  

\newcommand{\specsymb} {\sigma} 

\newcommand{\spec}[2][{}]   {\specsymb_{\mathrm{#1}}(#2)}

 \newcommand{\Err}{\mathcal O}

\def\XXint#1#2#3{{\setbox0=\hbox{$#1{#2#3}{\int}$}
     \vcenter{\hbox{$#2#3$}}\kern-.5\wd0}}

\newcommand{\R}{\mathbb{R}} 
\newcommand{\C}{\mathbb{C}} 
\newcommand{\N}{\mathbb{N}} 
\newcommand{\Sphere}{\mathbb{S}} 

\newcommand{\eps}{\varepsilon} 
\renewcommand{\phi}{\varphi}   
\newcommand{\e}{\mathrm e}  
\newcommand{\im}{\mathrm i} 

\newcommand{\wt}{\widetilde}           
\newcommand {\qf}[1]{\mathfrak{#1}}    

\newcommand{\HS}{\mathcal H}           

\newcommand{\Sobsymb} {\mathsf H}      
\newcommand{\Contsymb} {\mathsf C}     
\newcommand{\Lsymb}    {\mathsf L}     
\newcommand{\lsymb}    {\ell}          

\newcommand{\Cont}[2][{}]{\Contsymb^{#1}({#2})}

\newcommand{\Lsqr}[2][{}]{\Lsymb_2^{#1}({#2})} 
 
\newcommand{\lsqr}[2][{}]{\lsymb_2^{#1}({#2})}   

\newcommand{\Sob}[2][1]{\Sobsymb^{#1}({#2})}         

\newcommand{\Sobx}[3][1]{\Sobsymb_{{#2}}^{#1}({#3})} 
  
\newcommand{\abs}[1]{\lvert#1\rvert}   

\newcommand{\norm}[2][{}]{\|{#2}\|_{{#1}}}    
\newcommand{\normsqr}[2][{}]{\|{#2}\|^2_{#1}} 

\newcommand{\iprod}[3][{}]{\langle{#2},{#3}\rangle_{#1}}  

\newcommand{\set}[2]{\{ \, #1 \, | \, #2 \, \} }      
\newcommand{\bigset}[2]{\bigl\{ \, #1 \, \bigl|\bigr. \, #2 \, \bigr\} }

\newcommand{\map}[3]{ #1 \colon #2 \longrightarrow #3}    

\newcommand{\bd}  {\partial}                
\newcommand{\clo}[1]{\overline{{#1}}} 

\newcommand{\compl}[1]{#1^{\mathrm c}}       
\newcommand{\dcup}{\mathbin{\mathaccent\cdot\cup}}

\DeclareMathOperator*{\bigdcup}{\mathaccent\cdot{\bigcup}}

\newcommand{\restr}[1]{{\restriction}_{#1}} 

\newcommand{\conj}[1]{\overline {{#1}}}       
\newcommand{\orth}{\bot}                    

\newcommand{\1}{\mathbbm 1}                    
\newcommand{\und}{\quad\text{and}\quad}

\newcommand{\Neu}{{\mathrm N}}              
\newcommand{\Dir}{{\mathrm D}}              
\newcommand{\laplacian}[2][{}]{\Delta_{{#2}}^{{#1}}} 
\newcommand{\laplacianD}[1]{\laplacian[\Dir]{#1}} 
\newcommand{\laplacianN}[1]{\laplacian[\Neu]{#1}} 

\newcommand{\EW}[3][{}]{\lambda^{{#1}}_{#2}({#3})}
\newcommand{\EWN}[2]{\EW[\Neu]{#1}{#2}}      

\newcommand{\lapl} [1][{}]{\Delta_{#1}} 

\newcommand{\dlapl}[1][{}]{\pmb{\triangle}_{#1}}

\newcommand{\de} {\mathord{\mathrm d}} 

\newcommand{\dde}{\mathsf{d}}          
\newcommand{\dlaplacian}[2][{}]{\pmb{\triangle}_{{#2}}^{{#1}}}

\newcommand{\orient}[1]{\accentset{\curvearrowright}{#1}} 

\newcommand{\mc}{\mathcal}
\newcommand{\ul}{\underline}
\newcommand{\orul}[1]{\orient {\underline{#1}}}

\newcommand{\mbE}{\mathbb E} 

\newcommand{\Graph} X  
\newcommand{\Gmax}{\mc G^{\max}}

\newcommand{\stand}{\mathrm{std}}   

\newcommand{\vxeps}{{\eps,v}}
\newcommand{\edeps}{{\eps,e}}

\newcommand{\prim}{\mathrm{prim}}

\begin{document}
\title[Spectral analysis of metric graphs and related spaces]{Spectral
  analysis of metric graphs and related spaces}

\author{Olaf Post}      
\address{Institut f\"ur Mathematik,
         Humboldt-Universit\"at zu Berlin,
         Rudower Chaussee~25,
         12489 Berlin,
         Germany}
\email{post@math.hu-berlin.de}
\date{\today, \emph{File:} \jobname.tex}




\begin{abstract}
  The aim of the present article is to give an overview of spectral
  theory on metric graphs guided by spectral geometry on discrete
  graphs and manifolds.  We present the basic concept of metric graphs
  and natural Laplacians acting on it and explicitly allow infinite
  graphs. Motivated by the general form of a Laplacian on a metric
  graph, we define a new type of combinatorial Laplacian. With this
  generalised discrete Laplacian, it is possible to relate the
  spectral theory on discrete and metric graphs. Moreover, we describe
  a connection of metric graphs with manifolds. Finally, we comment on
  Cheeger's inequality and trace formulas for metric and discrete
  (generalised) Laplacians.
\end{abstract}

\maketitle

%
\section{Introduction}
\label{sec:intro}
%

A \emph{metric graph} $X$ is by definition a topological graph (i.e.,
a CW~complex of dimension $1$), where each edge $e$ is assigned a
length $\ell_e$. The resulting metric measure space allows to
introduce a family of ordinary differential operators acting on each
edge $e$ considered as interval $I_e=(0,\ell_e)$ with boundary
conditions at the vertices making the global operator self-adjoint.
One also refers to the pair of the graph and the self-adjoint
differential operator as \emph{quantum graph}.

Quantum graphs are playing an intermediate role between difference
operators on discrete graphs and partial differential operators on
manifolds. On the one hand, they are a good approximation of partial
differential operators on manifolds or open sets close to the graph,
see \Sec{mfd}. On the other hand, solving a system of ODEs reduces in
many cases to a discrete problem on the combinatorial graph, see
\Sec{rel.dg.mg}.

The spectral relation between metric and (generalised) discrete
Laplacians has the simplest form if the graph is \emph{equilateral},
i.e., if all lengths are the same, say, $\ell_e=1$. This fact and
related results have already been observed by many authors (see
e.g.~\cite{von-below:85,cattaneo:97,cartwright-woess:pre05,pankrashkin:06a,
  baker-faber:06,pankrashkin.talk:07,post:pre07c,bgp:08} and the
references therein).  Moreover, for non-equilateral graphs, one has at
least a spectral relation at the bottom of the spectrum. In
particular, one can define an \emph{index} (the Fredholm index of a
generalised ``exterior derivative'' in the discrete and metric case)
and show that they agree (\Thm{index}). The result extends the
well-known fact that the index equals the Euler characteristic for
standard graphs.  Such index formulas have been discussed e.g.\
in~\cite{kps:pre07,fkw:07,post:pre07a}. For convergence results of
a sequence of discrete Laplacians towards a metric graph Laplacian, we
refer to~\cite{faber:06} and the references therein.

Spectral graph theory is an active area of research.  We do not
attempt to give a complete overview here, and the choice of the
selected topics depends much on the author's taste. Results on
spectral theory of combinatorial Laplacians can be found e.g.\
in~\cite{dodziuk:84,mohar-woess:89,colin:98,cgy:96,chung:97,%
  higuchi-shirai:99,shirai:00,higuchi-shirai:04}.  For metric graph
Laplacians we mention the
works~\cite{roth:84,von-below:85,nicaise:87,kostrykin-schrader:99,%
  harmer:00b,kostrykin-schrader:03,kuchment:04,%
  friedman-tillich:pre04, kuchment:05,%
  baker-faber:06,kostrykin-schrader:06,pankrashkin:06a,%
  hislop-post:pre06,baker-rumely:07}.

Many concepts from spectral geometry on manifolds carry over to metric
and discrete graphs, and the right notion for a general scheme would
be a metric measure space with a Dirichlet form.  In particular,
metric graphs fall into this class; and they can serve as a toy model
in order to provide new results in spectral geometry.

This article is organised as follows: In the next section, we define
the generalised discrete Laplacians. \Sec{mg} is devoted to metric
graphs and their associated Laplacians. In \Sec{rel.dg.mg} we describe
relations between the discrete and metric graph Laplacians. \Sec{mfd}
contains the relation of a metric graph with a family of manifolds
converging to it. \Sec{cheeger} is devoted to the study of the first
non-zero eigenvalue of the Laplacian. In particular, we show Cheeger's
inequality for the (standard) metric graph Laplacian. \Sec{trace}
contains material on trace formulas for the heat operator associated
to (general) metric and discrete graph Laplacians. In particular, we
show a ``discrete path integral'' formula for generalised discrete
graph Laplacians (cf.~\Thm{trace.dg}).

\subsection*{Outlook and further developments}
Let us mention a few aspects which are not included in this article in
order to keep it at a reasonable size.  Our basic assumption is a
lower bound on the edge lengths. If we drop this condition, we obtain
\emph{fractal} metric graphs, i.e., (infinite) metric graphs with
$\inf_e \ell_e =0$. A simple example is given by a rooted tree, where
the length $\ell_n$ of an edge in generation $n$ tends to $0$. New
effects occur in this situation: for example, the Laplacian on
compactly supported functions can have more than one self-adjoint
extension; one needs additional boundary conditions at infinity (see
e.g.~\cite{solomyak:04}).

Another interesting subject are (infinite) covering graphs with finite
or compact quotient, for example Cayley graphs associated to a
finitely generated group. For example, if the covering group is
Abelian, one can reduce the spectral theory to a family of problems on
the quotient (with discrete spectrum) using the so-called
\emph{Floquet} theory. There are still open questions, for example
whether the (standard) discrete Laplacian of an equilateral maximal
Abelian covering has full spectrum or not. This statement is proven if
all vertices have even degree (using an ``Euler''-circuit). One can
ask whether similar statements hold also for general metric graph
Laplacians. For more details, we refer to~\cite{higuchi-shirai:99} and
the references therein.

Metric graphs have a further justification: The wave equation
associated to the (standard) metric graph Laplacian has finite
propagation speed, in contrast to the corresponding equation for the
(standard) discrete Laplacian (see~\cite[Sec.~4]{friedman-tillich:04}
for details). Note that the latter operator is bounded, whereas the
metric graph Laplacian is unbounded as differential operator.
Therefore, one can perform wave equation techniques on metric graphs
(and indeed, this has been done, see for example the scattering
approach in~\cite{kostrykin-schrader:99}).

\subsection*{Acknowledgements}
The author would like to thank the organisers of the programme
``Limits of graphs in group theory and computer science'' held at the
Bernoulli Center of the \emph{\'Ecole Polytechnique F\'ed\'erale de
  Lausanne} (EPFL), especially Prof.\ Alain Valette, for the kind
invitation and hospitality. The present article is an extended version
of a lecture held at the EPFL in March 2007.

%
\section{Discrete graphs and general Laplacians}
\label{sec:vx.sp}
%

In this section, we define a generalised discrete Laplacian, which
occurs also in the study of metric graph Laplacians as we will see in
\Sec{rel.dg.mg}.

\sloppy
Let us fix the notation: Suppose~$G$ is a countable, discrete,
weighted graph given by~$(V,E,\bd,\ell)$ where~$(V,E,\bd)$ is a usual
graph, i.e.,~$V$ denotes the set of vertices,~$E$ denotes the set of
edges,~$\map \bd E {V \times V}$ associates to each edge~$e$ the
pair~$(\bd_-e,\bd_+e)$ of its initial and terminal point (and
therefore an orientation). Abusing the notation, we also denote by
$\bd e$ the \emph{set} $\{\bd_-e,\bd_+e\}$.

  That~$G$ is an \emph{(edge-)weighted}
graph means that there is a \emph{length} or \emph{(inverse) edge
  weight function}~$\map \ell E {(0,\infty)}$ associating to each edge
$e$ a length~$\ell_e$. For simplicity, we consider \emph{internal}
edges only, i.e., edges of \emph{finite} length~$\ell_e < \infty$, and
we also make the following assumption on the lower bound of the edge
lengths:
\begin{assumption}
  \label{ass:len.bd}
  Throughout this article we assume that there is a constant $\ell_0 >
  0$ such that
  \begin{equation}
    \label{eq:len.bd}
    \ell_e \ge \ell_0, \qquad e \in E,
  \end{equation}
  i.e., that the weight function $\ell^{-1}$ is bounded.  Without loss
  of generality, we also assume that $\ell_0 \le 1$.
\end{assumption}
For each vertex~$v \in V$ we set
\begin{equation*}
  E_v^\pm := \set {e \in E} {\bd_\pm e = v} \qquad \text{and} \qquad
  E_v := E_v^+ \dcup E_v^-,
\end{equation*}
i.e.,~$E_v^\pm$ consists of all edges starting ($-$) resp.\ ending
($+$) at~$v$ and~$E_v$ their \emph{disjoint} union. Note that the
\emph{disjoint} union is necessary in order to allow self-loops, i.e.,
edges having the same initial and terminal point.  The \emph{degree
  of~$v \in V$} is defined as
\begin{equation*}
  \deg v := \card{E_v} = \card{E_v^+} + \card{E_v^-},
\end{equation*}
i.e., the number of adjacent edges at $v$. In order to avoid trivial
cases, we assume that~$\deg v \ge 1$, i.e., no vertex is isolated. We
also assume that $\deg v$ is finite for each vertex.

We want to introduce a vertex space allowing us to define Laplace-like
combinatorial operators motivated by general vertex boundary
conditions on quantum graphs. The usual discrete (weighted) Laplacian
is defined on \emph{scalar} functions $\map F V \C$ on the vertices
$V$, namely
\begin{equation}
  \label{eq:lap.std}
 \dlapl F(v) = - \frac 1 {\deg v} \sum_{e \in E_v} 
          \frac 1 {\ell_e} \bigl( F(v_e) - F(v) \bigr),
\end{equation}
where $v_e$ denotes the vertex on $e$ opposite to $v$. Note that
$\dlapl$ can be written as $\dlapl=\dde^* \dde$ with
\begin{equation}
  \label{eq:de.std}
  \map \dde {\lsqr V} {\lsqr E}, \qquad 
  (\dde F)_e = F(\bd_+ e) - F(\bd_- e),
\end{equation}
where $\lsqr V$ and $\lsqr E$ carry the norms defined by
\begin{equation}
  \label{eq:norm.std}
  \normsqr[\lsqr V] F
  := \sum_{v \in V} \abs{F(v)}^2 \deg v \und
  \normsqr[\lsqr E] \eta
  := \sum_{e \in E} \abs{\eta_e}^2\frac 1 {\ell_e},
\end{equation}
and $\dde^*$ denotes the adjoint with respect to the corresponding
inner products. We sometimes refer to functions in $\lsqr V$ and
$\lsqr E$ as \emph{$0$-} and \emph{$1$-forms}, respectively.

We would like to carry over the above concept for the vertex space
$\lsqr V$ to more general vertex spaces $\mc G$. The main motivation
to do so are metric graph Laplacians with general vertex boundary
conditions as defined in \Sec{mg} and their relations with discrete
graphs (cf.~\Sec{rel.dg.mg}).

\begin{definition}
  \label{def:vx.sp}
  \item
  \begin{enumerate}
  \item Denote by $\Gmax_v := \C^{E_v}$ the \emph{maximal vertex space
      at the vertex~$v \in V$}, i.e., a value~$\ul F(v) \in \Gmax_v$
    has~$\deg v$ components, one for each adjacent edge.  A (general)
    \emph{vertex space at the vertex $v$} is a linear subspace $\mc
    G_v$ of $\Gmax_v$.

  \item The corresponding (total) vertex spaces associated to the
    graph $(V,E,\bd)$ are
    \begin{equation*}
      \Gmax := \bigoplus_{v \in V} \Gmax_v \und
      \mc G := \bigoplus_{v \in V} \mc G_v,
    \end{equation*}
    respectively.  Elements of $\mc G$ are also called
    \emph{$0$-forms}.  The space $\mc G$ carries its natural Hilbert
    norm, namely
    \begin{equation*}
      \normsqr[\mc G] F
      := \sum_{v \in V} \abs{\ul F(v)}^2 
      = \sum_{v \in V} \sum_{e \in E_v}\abs{F_e(v)}^2.
    \end{equation*}
    Associated to a vertex space is an orthogonal projection~$P =
    \bigoplus_{v \in V} P_v$ in~$\Gmax$, where~$P_v$ is the orthogonal
    projection in~$\Gmax_v$ onto~$\mc G_v$.
  \item We call a general subspace $\mc G$ of $\Gmax$ \emph{local} iff
    it decomposes with respect to the maximal vertex spaces, i.e., if
    $\mc G = \bigoplus_v \mc G_v$ and $\mc G_v \le \Gmax_v$.
    Similarly, an operator $A$ on $\mc G$ is called \emph{local} if it
    is decomposable with respect to the above decomposition.
  \item The \emph{dual} vertex space associated to $\mc G$ is defined
    by~$\mc G^\orth := \Gmax \ominus \mc G$ and has
    projection~$P^\orth = \1 - P$.
  \end{enumerate}
\end{definition}

Note that a local subspace $\mc G$ is closed since $\mc G_v \le
\Gmax_v$ is finite dimensional.  Alternatively, a vertex space is
characterised by fixing an orthogonal projection~$P$ in~$\mc G$ which
is local.

\begin{example}
  \label{ex:vx.sp}
  The names of the vertex spaces in the examples below will become
  clear in the quantum graph case. For more general cases, e.g.\ the
  discrete magnetic Laplacian, we refer to~\cite{post:pre07a}.
  \begin{enumerate}
  \item 
    \label{cont}
    Choosing $\mc G_v = \C \ul \1(v)= \C(1, \dots, 1)$, we obtain the
    \emph{continuous} or \emph{standard} vertex space denoted by~$\mc
    G_v^\stand$.  The associated projection is
    \begin{equation*}
      P_v = \frac 1 {\deg v} \mbE
    \end{equation*}
    where~$\mbE$ denotes the square matrix of rank~$\deg v$ where all
    entries equal~$1$.  This case corresponds to the standard discrete
    case mentioned before. Namely, the natural identification
    \begin{equation*}
      \map{\wt \bullet}{\mc G^\stand}{\lsqr V}, \qquad
      F \mapsto \wt F, \qquad \wt F(v) := F_e(v),
    \end{equation*}
    (the latter value is independent of $e \in E_v$) is isometric,
    since the weighted norm in $\lsqr V$ and the norm in $\mc
    G^\stand$ agree, i.e., 
    \begin{equation*}
      \normsqr[\mc G^\stand] F 
      = \sum_{v \in V} \sum_{e \in E_v} \abs{F_e(v)}^2
      = \sum_{v \in V} \abs{\wt F(v)}^2 \deg v
      = \normsqr[\lsqr V] {\wt F}.
    \end{equation*}
    \item We call $\mc G_v^{\min} := 0$ the \emph{minimal} or
      \emph{Dirichlet} vertex space, similarly, $\Gmax$ is called the
      \emph{maximal} or \emph{Neumann} vertex space. The corresponding
      projections are $P=0$ and $P=\1$.
    \item
      \label{dim2}
      Assume that $\deg v=4$ and define a vertex space of dimension
      $2$ by
    \begin{equation*}
      \mc G_v = \C(1,1,1,1) \oplus \C(1,\im,-1,-\im).
    \end{equation*}
    The corresponding orthogonal projection is
    \begin{equation*}
      P = \frac 1 4
      \begin{pmatrix}
        2 & 1+\im & 0 & 1-\im\\
        1+\im & 0 & 1-\im & 2\\
        0 & 1-\im & 2 & 1+\im\\
        1-\im & 2 & 1+\im & 0
      \end{pmatrix}.
    \end{equation*}
    We will show some invariance properties of this vertex space in
    \Exenum{vx.sp.irr}{irr}.
  \end{enumerate}
\end{example}

For the next definition, we need some more notation. Let $E_{0,v}
\subset E_v$ be a subset of the set of adjacent edges at $v$.  We
denote by $\mc G_v \restr {E_{0,v}}$ the subspace of $\mc G_v$ where
the coordinates \emph{not} in $E_{0,v}$ are set to $0$, i.e.,
\begin{equation*}
  \mc G_v \restr {E_{0,v}} := 
        \set {\ul F(v)} 
             {F_e(v)=0, \quad \forall \; e \in E_v \setminus E_{0,v}}.
\end{equation*}
\begin{definition}
  A vertex space $\mc G_v$ at the vertex $v$ is called
  \emph{irreducible} if for any decomposition $E_v = E_{1,v} \dcup
  E_{2,v}$ such that $\mc G_v = \mc G_v \restr {E_{1,v}} \oplus \mc
  G_v \restr {E_{2,v}}$ we have either $E_{1,v}=\emptyset$ or
  $E_{2,v}=\emptyset$.  A vertex space $\mc G$ associated to a graph
  $G$ is \emph{irreducible} if all its components $\mc G_v$ are
  irreducible.
  
  By definition, the minimal vertex space $\mc G^{\min}_v=0$ is
  irreducible iff $\deg v=1$.
\end{definition}
In other words, a vertex space $\mc G_v$ is irreducible, if its
projection $P_v$ does not have block structure (in the given
coordinates).  The notion of \emph{irreducibility} is useful in order
to obtain a ``minimal'' representation of $\mc G$ by splitting a
vertex with a reducible vertex space into several vertices. Repeating
this procedure, we obtain:
\begin{lemma}
  \label{lem:dec}
  For any vertex space $\mc G$ associated to a graph $G=(V,E,\bd)$,
  there exists a graph $\wt G=(\wt V,E,\wt \bd)$ and a surjective
  graph morphism $\map \pi {\wt G} G$ such that $\mc G$ decomposes as
  \begin{equation*}
    \mc G = \bigoplus_{\wt v \in \wt V} \wt {\mc G}_{\wt v} \und
    \mc G_v = \bigoplus_{\wt v \in \pi^{-1}\{v\}} 
        \wt {\mc G}_{\wt v}.
  \end{equation*}
  In addition, each $\wt {\mc G}_{\wt v}$ is irreducible.
\end{lemma}
Note that the edge set of $\wt G$ is the same as for the original
graph $G$.
\begin{proof}
  We construct the vertex set $\wt V$ of $\wt G$ as follows: Let $v
  \in V$ and $\mc G_v$ be an irreducible vertex space, then $v$ is
  also an element of $\wt V$. Otherwise, if $\mc G_v = \mc G_v \restr
  {E_{1,v}} \oplus \mc G_v \restr {E_{2,v}}$ is a reducible vertex
  space (for $G$), we replace the vertex $v$ in $V$ by two different
  vertices $v_1$, $v_2$ in $\wt V$ with adjacent edges $E_{1,v}$ and
  $E_{2,v}$, in particular, $\wt G=(V \setminus \{v\} \cup
  \{v_1,v_2\}, E, \wt \bd)$ where
  \begin{equation*}
    \wt \bd_\pm e =
    \begin{cases}
      \bd_\pm e,& \text{if $\bd_\pm e \ne v$,}\\
      v_i, & \text{if $\bd_\pm = v$ and $e \in E_{v_i}$ for $i=1,2$.}
    \end{cases}
  \end{equation*}
  The associated vertex space at $v_i$ is $\wt {\mc G}_{v_i} := \mc
  G_v \restr {E_{i,v}}$ for $i=1,2$.  Repeating this procedure, we
  finally end with a graph $\wt G$ (denoted with the same symbol),
  such that each vertex space $\wt{\mc G}_{\wt v}$ is irreducible. The
  map $\pi$ is defined by $\pi e=e$ and $\pi \wt v = v$ if $\wt v$ came from
  splitting a vertex space at the original vertex $v$. It is easy to
  see that $\pi$ is a graph morphism (i.e, $\bd_\pm \pi e = \pi \wt
  \bd _\pm e$) and surjective.
\end{proof}
\begin{definition}
  \label{def:irr}
  We call the graph $\wt G$ constructed in \Lem{dec} the
  \emph{irreducible} graph of the vertex space $\mc G$ associated to
  the graph $G$. We say that the vertex space is \emph{connected} if
  the associated irreducible graph is a connected graph.
\end{definition}
Note that on the level of the vertex space $\mc G$, passing to the
irreducible graph is just a reordering of the coordinate labels,
namely, a regrouping of the labels into smaller sets.

For example, the maximal vertex space $\Gmax$ associated to a graph
$G$ (with $\deg v \ge 2$ for all vertices $v$) is \emph{not}
irreducible, and its irreducible graph is
\begin{equation}
  \label{eq:dec.gr}
  \wt G = \bigdcup_{e \in E} G_e \qquad
  \text{where} \qquad G_e:=(\bd e,\{e\}, \bd \restr {\{e\}})
\end{equation}
is a graph with two vertices and one edge only. The vertex space is
\begin{equation}
  \label{eq:g.max}
  \Gmax = \bigoplus_{e \in E} (\C_{\bd_-e} \oplus \C_{\bd_+e})
\end{equation}
where $\C_{\bd_\pm e}$ is a copy of $\C$.  The irreducible graph of
the minimal vertex space $\mc G^{\min}=0$ is the same as above.

However, the standard vertex space $\mc G^\stand$ associated to a
graph $G$ is already irreducible and $\wt G=G$. Therefore, the
standard vertex space is connected iff the underlying graph is
connected; i.e., the notion of ``connectedness'' agrees with the usual
one.


Now, we define a generalised \emph{coboundary operator} or
\emph{exterior derivative} associated to a vertex space. We use this
exterior derivative for the definition of an associated Laplace
operator below:

\begin{definition}
  \label{def:discr.ext.der}
  Let~$\mc G$ be a vertex space of the graph~$G$. The
  \emph{exterior derivative} on~$\mc G$ is defined via
  \begin{equation*}
    \map{\dde_{\mc G}}{\mc G}
              {\lsqr E}, \qquad
    (\dde_{\mc G}  F)_e := F_e(\bd_+ e) - F_e(\bd_- e),
  \end{equation*}
  mapping~$0$-forms onto~$1$-forms.
\end{definition}
We often drop the subscript~$\mc G$ for the vertex space.  The proof
of the next lemma is straightforward (see
e.g.~~\cite[Lem.~3.3]{post:pre07a}):
\begin{lemma}
  \label{lem:discr.ext.der}
  Assume the lower lengths bound~\eqref{eq:len.bd}, then~$\dde$ is
  norm-bounded by~$\sqrt {2/\ell_0}$. The adjoint
  \begin{equation*}
    \map {\dde^*}{\lsqr E}{\mc G}
  \end{equation*}
  fulfills the same norm bound and is given by
  \begin{equation*}
    (\dde^* \eta)(v) 
    = P_v \Bigl( \Bigl\{ \frac 1 \ell_e \orient \eta_e(v) \Bigr\} \Bigr)
    \in \mc G_v,
  \end{equation*}
  where~$\orient \eta_e(v):= \pm \eta_e$ if~$v=\bd_\pm e$ denotes the
  \emph{oriented} evaluation of~$\eta_e$ at the vertex~$v$.
\end{lemma}

\begin{definition}
  \label{def:discr.laplace}
  The \emph{discrete generalised Laplacian} associated to a vertex
  space $\mc G$ is defined as $\dlaplacian {\mc G} := \dde_{\mc G}^*
  \dde_{\mc G}$, i.e.,
  \begin{equation*}
    (\dlaplacian {\mc G} F)(v)
    = P_v \Bigl( \Bigl\{ \frac 1 \ell_e
         \bigl( F_e(v) - F_e(v_e) \bigr) \Bigr\} \Bigr)
  \end{equation*}
  for $F \in \mc G$, where~$v_e$ denotes the vertex on $e \in E_v$
  opposite to~$v$.
\end{definition}

\begin{remark}
  \indent
  \begin{enumerate}
  \item From \Lem{discr.ext.der} it follows that $\dlaplacian {\mc G}$
    is a bounded operator on $\mc G$ with norm estimated from above
    by~$2/\ell_0$.

  \item Note that the orientation of the edges plays no role for the
    ``second order'' operator $\dlaplacian {\mc G}$.

  \item We can also define a Laplacian $\dlaplacian[1] {\mc G}:=
    \dde_{\mc G} \dde_{\mc G}^*$ acting on the space of ``$1$-forms''
    $\lsqr E$ (and $\dlaplacian[0] {\mc G}:= \dlaplacian {\mc G} =
    \dde_{\mc G}^* \dde_{\mc G}$).  For more details and the related
    supersymmetric setting, we refer to~\cite{post:pre07a}.  
    In particular, we have
    \begin{equation*}
       \spec {\dlaplacian[1] {\mc G}} \setminus \{0\}
        = \spec {\dlaplacian[0]{\mc G}} \setminus \{0\}.
    \end{equation*}
    Moreover, in~\cite[Ex.~3.16--3.17]{post:pre07a} we discussed how
    these generalised Laplacians can be used in order to analyse the
    (standard) Laplacian on the line graph and subdivision graph
    associated to $G$ (see also~\cite{shirai:00}).
  \item Assume that $G$ is equilateral (i.e., $\ell_e=1$), which
    implies $\spec{\laplacian{\mc G}}\subseteq [0,2]$. Then using the
    $1$-form Laplacian, one can show the spectral relation
    \begin{equation*}
      \spec {\dlaplacian {\mc G^\orth}} \setminus \{0,2\}
      = 2 -(\spec {\dlaplacian{\mc G}} \setminus \{0,2\}),
    \end{equation*}
    i.e., if~$\lambda \notin \{0,2 \}$, then~$\lambda \in \spec
    {\dlaplacian{\mc G^\orth}}$ iff~$2-\lambda \in \spec
    {\dlaplacian{\mc G}}$ (cf.~\cite[Lem.~3.13~(iii)]{post:pre07a}).
  \end{enumerate}
\end{remark}
The next example shows that we have indeed a generalisation of the
standard discrete Laplacian:
\begin{example}
  \label{ex:ext.der}
  \indent
  \begin{enumerate}
  \item For the standard vertex space~$\mc G^\stand$, it is convenient
    to use the unitary transformation from $\mc G^\stand$ onto $\lsqr
    V$ associating to $F \in \mc G$ the (common value) $\wt F(v):=
    F_e(v)$ as in \Exenum{vx.sp}{cont}.  Then the exterior derivative
    and its adjoint are unitarily equivalent to
    \begin{equation*}
      \map{\wt \dde}{\lsqr V}{\lsqr E}, \qquad
      (\wt \dde \wt F)_e = \wt F(\bd_+ e) - \wt F(\bd_- e)
    \end{equation*}
    and
    \begin{equation*}
      (\wt \dde^* \eta)(v) 
      = \frac 1 {\deg v} \sum_{e \in E_v} 
      \frac 1 {\ell_e} \orient \eta_e(v),
    \end{equation*}
    i.e.,~$\wt \dde$ is the classical coboundary operator already
    defined in~\eqref{eq:de.std} and~$\wt \dde^*$ its adjoint.

    Moreover, the corresponding discrete Laplacian $\dlaplacian{\mc
      G^\stand}$ is unitarily equivalent to the usual discrete
    Laplacian~$\dlapl=\wt \dde^* \wt \dde$ defined
    in~\eqref{eq:lap.std} as one can easily check.

  \item Passing to the irreducible graph of a vertex space $\mc G$ is
    a reordering of the coordinate labels, and in particular, the
    Laplacian is the same (up to the order of the coordinate labels).
    Namely, for the minimal vertex space $\mc G^{\min}=0$, we
    have~$\dde=0$, $\dde^*=0$ and $\dlaplacian{\mc G^{\min}}=0$.

    For the maximal vertex space, we have
    \begin{equation*}
        (\dlaplacian{\mc G^{\max}} F)_e(v) = 
          \Bigl\{ 
               \frac 1 {\ell_e} \bigl(F_e(v) - F_e(v_e) \bigr) 
          \Bigr\}_{e \in E_v}
    \end{equation*}
    and
    \begin{equation*}
        \dlaplacian{\mc G^{\max}} = 
        \bigoplus_{e \in E} \dlaplacian{G_e} \qquad
        \text{where} \qquad
        \dlaplacian{G_e} \cong \frac 1 {\ell_e}
        \begin{pmatrix}
          1 & -1 \\ -1 & 1
        \end{pmatrix}.
    \end{equation*}
    In particular, in both cases, the Laplacians are decoupled and any
    connection information of the graph is lost.
  \end{enumerate}
\end{example}
Of course, the decoupled minimal and maximal cases are uninteresting
when analysing the graph and its properties. Moreover, it is natural
to assume that the vertex space is connected and irreducible, since
the other cases can be reduced to this one.

Let us analyse the generalised Laplacian in the special case when all
lengths are equal, say, $\ell_e=1$ and when there are no double edges.
Then we can write the Laplacian in the form
\begin{equation*}
  \laplacian {\mc G} = \1 - M_{\mc G}, \qquad M_{\mc G}:=P A^{\max},
\end{equation*}
where $\map{M_{\mc G}}{\mc G}{\mc G}$ is called the \emph{principle
  part} of the generalised discrete Laplacian, and $\map {A^{\max}}
\Gmax \Gmax$ the \emph{generalised adjacency matrix}, defined by
\begin{equation*}
  A^{\max} \{F(w)\}_w = \{A^{\max}(v,w) F(w)\}_v, 
    \quad
  \map{A^{\max}(v,w)}{\C^{E_w}} {\C^{E_v}}
\end{equation*}
for $F \in \Gmax$. Furthermore, $A^{\max}(v,w)=0$ if $v,w$ are not joined
by an edge and
\begin{equation*}
  A^{\max}(v,w)_{e,e'}=\delta_{e,e'}, \quad e \in E_v,\; e' \in E_w
\end{equation*}
otherwise. In particular, written as a matrix, $A^{\max}(v,w)$ has
only one entry $1$ and all others equal to $0$. The principle part of
the Laplacian then has the form
\begin{equation*}
  (M_{\mc G} F)(v) = 
  \sum_{e \in E_v} A_{\mc G}(v,v_e) F(v_e) ,
\end{equation*}
for $F \in \mc G$ similar to the form of the principle part of the
standard Laplacian defined for $\mc G^\stand \cong \lsqr V$, where
\begin{equation*}
  \map{A_{\mc G}(v,w):= P_v A^{\max}(v,w) P_w}{\mc G_w}{\mc G_v}.
\end{equation*}
Equivalently,
\begin{equation}
  \label{eq:markov}
  M_{\mc G} = 
     \bigoplus_{v \in V} \sum_{w \in V} A_{\mc G}(v,w) 
\end{equation}
where the sum is actually only over those vertices $w$ connected with
$v$.  In particular, in the standard case $\mc G=\mc G^\stand$, the
matrix $A_{\mc G^\stand}(v,w)$ consists of one entry only since $\mc
G^\stand_v \cong \C (\deg v)$ isometrically, namely $A_{\mc
  G^\stand}(v,w)=1$ if $v$ and $w$ are connected and $0$ otherwise,
i.e., $A_{\mc G^\stand}$ is (unitarily equivalent to) the standard
adjacency operator in $\lsqr V$.

Let us return to the general situation (i.e., general lengths $\ell_e$
and possibly double edges).  In~\cite[Lem.~2.13]{post:pre07a} we
showed the following result on symmetry of a vertex space:
\begin{lemma}
  \label{lem:inv.sym}
  Assume that the vertex space $\mc G_v$ of a vertex $v$ with degree
  $d=\deg v$ is invariant under permutations of the coordinates $e \in
  E_v$, then $\mc G_v$ is one of the spaces $\mc G_v^{\min}=0$, $\mc
  G^{\max} =\C^{E_v}$, $\mc G^\stand=\C(1,\dots,1)$ or $(\mc
  G^\stand)^\orth$, i.e., only the minimal, maximal, standard and dual
  standard vertex space are invariant.
\end{lemma}

If we only require invariance under the cyclic group of order $d$, we
have the following result:
\begin{lemma}
  \label{lem:inv.cyc}
  Assume that the vertex space $\mc G_v$ of a vertex $v$ with degree
  $d=\deg v$ is invariant under a cyclic permutation of the
  coordinates $e \in E_v=\{e_1,\dots,e_d\}$, i.e., edge $e_i \mapsto
  e_{i+1}$ and $e_d \mapsto e_1$, then $\mc G_v$ is an orthogonal sum
  of spaces of the form $\mc G_v^p = \C(1,\theta^p, \theta^{2p},
  \dots, \theta^{(d-1)p})$ for $p=0, \dots, d-1$, where $\theta =
  \e^{2 \pi \im/d}$.
\end{lemma}
\begin{proof}
  The (representation-theoretic) irreducible vector spaces invariant
  under the cyclic group are one-dimensional (since the cyclic group
  is Abelian) and have the form $\mc G_v^p$ as given below.
\end{proof}
We call $\mc G_v^p$ a \emph{magnetic} perturbation of $\mc
G_v^\stand$, i.e., the components of the generating vector $(1, \dots,
1)$ are multiplied with a phase factor (see
e.g.~\cite[Ex.~2.10~(vii)]{post:pre07a}).
\begin{example}
  \label{ex:vx.sp.irr}
\indent
  \begin{enumerate}
  \item If we require that the vertex space $\mc G_v$ is cyclic
    invariant with \emph{real} coefficients in the corresponding
    projections, then $\mc G_v$ is $\C(1,\dots,1)$ or
    $\C(1,-1,\dots,1,-1)$ (if $d$ even) or their sum.  But the sum is
    reducible since
    \begin{equation*}
      \qquad\qquad 
      \mc G_v = \C(1,\dots,1) \oplus \C(1,-1,\dots,1,-1) =
      \C(1,0,1,0,\dots,1,0) \oplus \C(0,1,0,1,\dots,0,1)
    \end{equation*}
    and the latter two spaces are standard with degree $d/2$. In other
    words, the irreducible graph at $v$ associated to the boundary
    space $\mc G_v$ splits the vertex $v$ into two vertices $v_1$ and
    $v_2$ adjacent with the edges with even and odd labels,
    respectively. The corresponding vertex spaces are standard.
  \item
    \label{irr}
    The sum of two cyclic invariant spaces is not always
    reducible: Take the cyclic invariant vertex space $\mc G_v =\mc
    G_v^0 \oplus \mc G_v^1 \le \C^4$ of dimension $2$ given in
    \Exenum{vx.sp}{dim2}.  Note that $\mc G_v$ is irreducible, since
    the associated projection $P$ does not have block structure.  This
    vertex space is maybe the simplest example of an (cyclic
    invariant) irreducible vertex space which is not standard or dual
    standard. Note that if $\deg v=3$, then an irreducible vertex
    space is either standard or dual standard (or the corresponding
    version with ``weights'', i.e., $(1,\dots,1)$ replaced by a
    sequence of non-zero numbers).
  \end{enumerate}
\end{example}

We finally develop an index theory associated to a vertex space $\mc
G$.  We define the Hilbert chain associated to a vertex space $\mc G$
as
\begin{equation*}
  \mc C_{G,\mc G} \colon 0 
   \longrightarrow \mc G \stackrel {\dde_\mc G} 
   \longrightarrow  \lsqr E \longrightarrow 0.
\end{equation*}
Obviously, the chain condition is trivially satisfied since only one
operator is non-zero. In this situation and since we deal with Hilbert
spaces, the associated cohomology spaces (with coefficients in $\C$)
can be defined as
    \begin{align*}
      H^0(G,\mc G) &:= \ker \dde_{\mc G} \cong \ker \dde_{\mc G}/\ran 0, \\
      H^1(G,\mc G) &:= \ker \dde^*_{\mc G} = \ran \dde_{\mc G}^\orth
      \cong \ker 0/\ran \dde_{\mc G}
    \end{align*}
where $\ran A:=A(\HS_1)$ denotes the range (``image'') of the operator
$\map A {\HS_1}{\HS_2}$. The \emph{index} or \emph{Euler
  characteristic} of this cohomology is then defined as
\begin{equation*}
  \ind (G,\mc G) := \dim \ker \dde_{\mc G} - \dim \ker \dde_{\mc G}^*,
\end{equation*}
i.e., the \emph{Fredholm index} of $\dde_{\mc G}$, provided at least
one of the dimensions is finite.  Note that for the standard vertex
space $\mc G^\stand \cong \lsqr V$, the exterior derivative is just
(equivalent to) the classical coboundary operator defined
in~\eqref{eq:de.std}. In particular, the corresponding homology spaces
are the classical ones, and $\dim H^p(G,\mc G^\stand)$ counts the
number of components ($p=0$) and edges not in a spanning tree ($p=1$).

Using the stability of the index under (at least) continuous
perturbations, we can calculate the index via simple (decoupled) model
spaces and obtain (see~\cite[Sec.~4]{post:pre07a}):
\begin{theorem}
  \label{thm:ind.dg}
  Let $\mc G$ be a vertex space associated with the finite graph
  $G=(V,E,\bd)$, then
  \begin{equation*}
    \ind (G,\mc G) = \dim \mc G - \card E.
  \end{equation*}
\end{theorem}
Note that in particular, if $\mc G=\mc G^\stand$, i.e., if $\mc G
\cong \lsqr V$ is the standard vertex space, we recover the well-known
formula for (standard) discrete graphs, namely
\begin{equation*}
  \ind (G,\mc G^\stand) = \card V - \card E,
\end{equation*}
i.e., the index is the Euler characteristic $\chi(G):=\card V - \card
E$ of the graph $G$.  On the other hand, in the ``extreme'' cases, we
have
\begin{equation*}
  \ind (G,\Gmax) = \card E 
               \und
  \ind (G,\mc G^{\min}) = - \card E.
\end{equation*}
since $\dim \Gmax = \sum_{v \in V} \deg v = 2\card E$ and $\dim \mc
G^{\min}=0$. Again, the index equals the Euler characteristic of the
decoupled graph $\chi(\bigdcup_e G_e))=\sum_e \chi(G_e)=2\card E$ (see
\Eq{dec.gr}) resp.\ the relative Euler characteristic
$\chi(G,V)=\chi(G)-\chi(V)=-\card E$.

In \cite[Lem.~4.4]{post:pre07a} we established a general result on the
cohomology of the dual~$\mc G^\orth$ of a vertex space~$\mc G$. It
shows that actually, $\mc G^\orth$ and the \emph{oriented} version of
$\mc G$, i.e., $\orient {\mc G}=\set{F \in \Gmax}{\orient F \in \mc
  G}$, are related:
\begin{lemma}
  \label{lem:ind.dual}
  Assume that the global length bound
  \begin{equation}
  \label{eq:len.2bd}
    \ell_0 \le \ell_e \le \ell_+ \qquad \text{for all $\e \in E$}
  \end{equation}
  holds for some constants $0 < \ell_0 \le \ell_+ < \infty$. Then
  $H^0(G,\mc G^\orth)$ and $H^1(G, \orient{\mc G})$
  are isomorphic.  In particular, if $G$ is finite, then $\ind(G,\mc
  G^\orth)= - \ind(G,\orient{\mc G})$.
\end{lemma}
The change of orientation also occurs in the metric graph case, see
e.g.~\Lem{adj}.

\section{Metric graphs}
\label{sec:mg}
In this section, we fix the basic notion for metric and quantum graphs
and derive some general assertion needed later on.

\begin{definition}
  \label{def:mg}
  Let $G=(V,E,\bd)$ be a discrete graph. A \emph{topological graph}
  associated to $G$ is a CW complex $X$ containing only $0$-cells and
  $1$-cells, such that the $0$-cells are the vertices $V$ and the
  $1$-cells are labelled by the edge set $E$.

  A \emph{metric graph} $X=X(G,\ell)$ associated to a weighted
  discrete graph $(V,E,\bd,\ell)$ is a topological graph associated to
  $(V,E,\bd)$ such that for every edge $e \in E$ there is a continuous
  map $\map{\Phi_e}{\clo I_e}X$, $I_e:=(0,\ell_e)$, whose image is the
  $1$-cell corresponding to $e$, and the restriction
  $\map{\Phi_e}{I_e}{\Phi(I_e)} \subset X$ is a homeomorphism. The
  maps $\Phi_e$ induce a metric on $X$. In this way, $X$ becomes a
  metric space.
\end{definition}
Given a weighted discrete graph, we can abstractly construct the
associated metric graph as the disjoint union of the intervals $I_e$
for all $e \in E$ and appropriate identifications of the end-points of
these intervals (according to the combinatorial structure of the
graph), namely
\begin{equation}
  \label{eq:constr.mg}
  X = \bigdcup_{e \in E} \clo I_e / {\sim}.
\end{equation}
We denote the union of the $0$-cells and the union of the (open)
$1$-cells (edges) by $X^0$ and $X^1$, i.e.,
\begin{equation*}
  X^0 = V \hookrightarrow X, \qquad 
  X^1 =\bigcup_{e \in E} I_e \hookrightarrow X,
\end{equation*}
and both subspaces are canonically embedded in $X$.
\begin{remark}
  \label{rem:met.space}
  \indent
  \begin{enumerate}
  \item 
    \label{coord}
    The metric graph $X$ becomes canonically a \emph{metric measure
      space} by defining the distance of two points to be the length
    of the shortest path in $X$, joining these points. We can think of
    the maps $\map{\Phi_e}{I_e} X$ as coordinate maps and the Lebesgue
    measures on the intervals $I_e$ induce a (Lebesgue) measure on the
    space $X$. We will often abuse the notion and write $X=(G,\ell)$
    or $X=(V,E,\bd,\ell)$ for the metric graph associated to the
    weighted discrete graph $(G,\ell)$ with $G=(V,E,\bd)$.
  \item Note that two metric graphs $X=(G,\ell)$, $X'=(G',\ell')$ can
    be isometric as metric spaces but not isomorphic as graphs: The
    metric on a metric graph $X$ cannot distinguish between a single
    edge $e$ of length $\ell_e$ in $G$ and two edges $e_1$, $e_2$ of
    length $\ell_{e_1}$ , $\ell_{e_2}$ with
    $\ell_e=\ell_{e_1}+\ell_{e_2}$ joined by a single vertex of degree
    $2$ in $G'$: The underlying graphs are not (necessarily)
    isomorphic.  For a discussion on this point, see for
    example~\cite[Sec.~2]{baker-rumely:07}.
  \end{enumerate}
\end{remark}

Since a metric graph is a topological space, and isometric to
intervals outside the vertices, we can introduce the notion of
measurability and differentiate function on the edges.  We start with
the basic Hilbert space
\begin{gather*}
  \Lsqr X := \bigoplus_{e \in E} \Lsqr {I_e}, \qquad
  f = \{f_e\}_e \quad \text{with $f_e \in \Lsqr{I_e}$ and}\\
  \normsqr f = \normsqr[\Lsqr X] f 
  := \sum_{e \in E} \int_{I_e} \abs{f_e(x)}^2 \dd x.
\end{gather*}

In order to define a natural Laplacian on $\Lsqr X$ we introduce
the \emph{maximal} or \emph{decoupled} Sobolev space of order $k$ as
\begin{gather*}
  \Sobx[k] \max X
    := \bigoplus_{e \in E} \Sob [k]{I_e},\\
  \normsqr[{\Sobx[k] \max X}] f 
  := \sum_{e \in E} \normsqr[{\Sob [k]{I_e}}] {f_e},
\end{gather*}
where $\Sob[k] {I_e}$ is the classical Sobolev space on the interval
$I_e$, i.e., the space of functions with (weak) derivatives in $\Lsqr
{I_e}$ up to order $k$.  We define the \emph{unoriented} and
\emph{oriented} value of $f$ on the edge $e$ at the vertex $v$ by
\begin{equation*}
  \ul f_e(v) :=
  \begin{cases}
    f_e(0), & \text{if $v=\bd_-e$,}\\
    f_e(\ell(e)), & \text{if $v=\bd_+e$}
  \end{cases} \und
  \orul f_e(v) :=
  \begin{cases}
    -f_e(0), & \text{if $v=\bd_-e$,}\\
    f_e(\ell(e)), & \text{if $v=\bd_+e$}.
  \end{cases}
\end{equation*}
Note that $\ul f_e(v)$ and $\orul f_e(v)$ are defined for $f \in \Sobx
\max X$. Even more, we have shown in~\cite[Lem.~5.2]{post:pre07a} the
following result:
\begin{lemma}
  \label{lem:mg.bd.eval}
  Assume the lower lengths bound~\eqref{eq:len.bd}, then the
  evaluation operators
  \begin{gather*}
    \map{\ul \bullet}{\Sobx \max X} \Gmax 
            \und 
    \map{\orul \bullet}{\Sobx \max X} \Gmax,
  \end{gather*}
  given by $f \mapsto \ul f=\{\{\ul f_e(v)\}_{e \in E_v}\}_v \in
  \Gmax=\bigoplus_v \Gmax_v = \bigoplus_v \C^{E_v}$ and similarly
  $\orul f \in \Gmax$, are bounded by $2 \ell_0^{-1/2}$.
\end{lemma}
These two evaluation maps allow a very simple formula of a partial
integration formula on the metric graph, namely
\begin{equation}
  \label{eq:part.int}
  \iprod[\Lsqr X] {f'} g = \iprod[\Lsqr X] {f}{-g'} + 
          \iprod[\Gmax] {\ul f}{\orul g},
\end{equation}
where $f'=\{f_e'\}_e$ and similarly for $g$.  Basically, this follows
from partial integration on each interval $I_e$ and a reordering of
the labels by
\begin{equation*}
  E = \bigdcup_{v \in V} E_v^+ = \bigdcup_{v \in V} E_v^-.
\end{equation*}
\begin{remark}
  If we distinguish between functions ($0$-forms) and vector fields
  ($1$-forms), we can say that $0$-forms are evaluated
  \emph{unoriented}, whereas $1$-forms are evaluated \emph{oriented}.
  In this way, we should interprete $f'$ and $g$ as $1$-forms and $f$,
  $g'$ as $0$-forms.
\end{remark}

Let $\mc G$ be a vertex space (i.e., a local subspace of $\Gmax$, or
more generally, a closed subspace) associated to the underlying
discrete graph. We define
\begin{equation*}
  \Sobx[k] {\mc G} X := 
  \bigset {f \in \Sobx[k] \max X} {\ul f \in \mc G} 
     \und
  \Sobx[k] {\orient{\mc G}} X := 
  \bigset {f \in \Sobx[k] \max X} {\orul f \in \mc G} .
\end{equation*}
Note that these spaces are closed in $\Sobx[k] \max X$ as
pre-image of the bounded operators $\ul \bullet$ and $\orul \bullet$,
respectively, of the closed subspace $\mc G$, and therefore itself
Hilbert spaces.

We can now mimic the concept of exterior derivative:
\begin{definition}
  \label{def:ext.der}
  The \emph{exterior derivative} associated to a metric graph $X$ and
  a vertex space $\mc G$ is the unbounded operator $\de_{\mc G}$ in
  $\Lsqr X$ defined by $\de_{\mc G} f:= f'$ for $f \in \dom \de_{\mc
    G} :=\Sobx {\mc G} X$.
\end{definition}
\begin{remark}
  \label{rem:ext.der}
  \indent
  \begin{enumerate}
  \item
    \label{closed}
    Note that $\de_{\mc G}$ is a closed operator (i.e., its graph is
    closed in $\Lsqr X \oplus \Lsqr X$), since $\Sobx {\mc G} X$ is a
    Hilbert space and the graph norm of $\de=\de_{\mc G}$ given by
    $\normsqr[\de] f := \normsqr {\de f} + \normsqr f$ is the Sobolev
    norm, i.e, $\norm [\de] f = \norm[\Sobx \max X] f$.
  \item We can think of $\de$ as an operator mapping $0$-forms into
    $1$-forms. Obviously, on a one-dimensional \emph{smooth} space,
    there is no need for this distinction, but the distinction between
    $0$- and $1$-forms makes sense through the boundary conditions
    $\ul f \in \mc G$, see also the next lemma.
  \end{enumerate}
\end{remark}
The adjoint of $\de_{\mc G}$ can easily be calculated from the partial
integration formula~\eqref{eq:part.int}, namely the boundary term has
to vanish for functions in the domain of $\de_{\mc G}^*$:
\begin{lemma}
  \label{lem:adj}
  The adjoint of $\de_{\mc G}$ is given by $\de_{\mc G}^* g = -g'$
  with domain $\dom \de_{\mc G}^* = \Sobx {\orient {\mc G}^\orth} X$.
\end{lemma}

As for the discrete operators, we define the Laplacian as
\begin{equation*}
  \laplacian{\mc G} := \de_{\mc G}^*\de_{\mc G}
\end{equation*}
with domain $\dom \laplacian{\mc G}:=\set{f \in \dom \de_{\mc G}} {\de
  f \in \dom \de_{\mc G}^*}$. Moreover, we have (see
e.g.~\cite[Thm.~17]{kuchment:04} or~\cite[Sec.~5]{post:pre07c} for
different proofs):
\begin{proposition}
  \label{prp:sa}
  Assume the lower lengths bound~\eqref{eq:len.bd},
  then~$\laplacian{\mc G}$ is self-adjoint on
  \begin{equation*}
    \dom \laplacian{\mc G}
       := \bigset{f \in \Sobx[2] \max X } 
                 {\ul f \in \mc G, \quad \orul f' \in \mc G^\orth}.
  \end{equation*}
\end{proposition}
\begin{proof}
  By definition of $\laplacian{\mc G}$, the Laplacian is the
  non-negative operator associated to the non-negative quadratic form
  $f \mapsto \normsqr {\de f}$ with domain $\Sobx {\mc G} X$. The
  latter is closed since $\Sobx {\mc G} X$ is a Hilbert space equipped
  with the associated quadratic form norm defined by $\normsqr[\Sob X]
  f = \normsqr{\de f} + \normsqr f$, see \Remenum{ext.der}{closed}.
  It remains to show that $\laplacian{\mc G}$ is a closed operator,
  i.e., $\dom \laplacian{\mc G}$ is a Hilbert space equipped with the
  graph norm defined by $\normsqr[\lapl] f := \normsqr f +
  \normsqr{f''}$.  By \Lem{mg.bd.eval}, the domain is a closed
  subspace of $\Sobx [2] \max X$, and it remains to show that the
  Sobolev and the graph norms
  \begin{equation*}
    \normsqr[{\Sobx[2]\max
    X}] f = \normsqr f + \normsqr {f'} +
  \normsqr {f''} \und \normsqr[\lapl]
  f = \normsqr f + \normsqr{f''},
  \end{equation*}
  are equivalent, i.e., that there is a constant $C>0$ such that
  $\normsqr {f'} \le C (\normsqr f + \normsqr {f''})$. The latter
  estimate is true under the global lower bound on the length
  function~\eqref{eq:len.bd} (see
  e.g.~\cite[App.~C]{hislop-post:pre06}).
\end{proof}

\begin{definition}
  \label{def:qg} A metric graph $X$ together with a self-adjoint
  Laplacian (i.e, an operator acting as $(\lapl f)_e = -f''_e$ on each
  edge) will be called \emph{quantum graph}.
\end{definition}
For example, $(X,\laplacian{\mc G})$ is a quantum graph; defined by
the data $(V,E,\bd, \ell, \mc G)$.
\begin{example}
  \label{ex:stand}
  The standard vertex space $\mc G^\stand$ leads to \emph{continuous}
  functions in $\Sobx {\mc G^\stand} X$, i.e., the value of $\ul
  f_e(v)$ is \emph{independent} of $e \in E_v$. Note that on each
  edge, we already have the embedding $\Sob {I_e} \subset \Cont
  {I_e}$, i.e., $f$ is already continuous inside each edge. In
  particular, a function $f$ is in the domain of $\laplacian{\mc
    G^\stand}$ iff $f \in \Sobx[2] \max X$, $f$ is continuous and
  $\orul f'(v) \in (\mc G^\stand)^\orth$. The latter condition on the
  derivative is a \emph{flux} condition, namely
  \begin{equation*}
    \sum_{e \in E_v} \orul f'_e(v) = 0
  \end{equation*}
  for all $v \in V$. The corresponding metric graph Laplacian
  $\laplacian {\mc G^\stand}$ is called \emph{standard}, or sometimes
  also \emph{Kirchhoff} Laplacian.
\end{example}
\begin{remark}
  \label{rem:sa}
  \indent
  \begin{enumerate}
  \item
    \label{other.sa}
    There are other possibilities how to define self-adjoint
    extensions of a Laplacian, namely for any self-adjoint (bounded)
    operator $L$ on $\mc G$, one can show that $\laplacian{(\mc G,L)}$
    is self-adjoint on
    \begin{equation*}
      \dom \laplacian{(\mc G,L)} :=
      \bigset {f \in \Sobx[2]{\mc G} X} {P\orul f' = L \ul f},
    \end{equation*}
    where $P$ is the projection in $\Gmax$ onto the space $\mc G$.
    The domain mentioned in \Prp{sa} corresponds to the case $L=0$.
    For more details, we refer e.g.\ to~\cite[Thm.~17]{kuchment:04}
    or~\cite[Sec.~4]{post:pre07c}, \cite{kps:pre07} (and references
    therein) and the next remark for another way of a parametrisation
    of self-adjoint extensions.
  \item 
    \label{scatt}
    One can encode the vertex boundary conditions also in a (unitary)
    operator $S$ on $\Gmax$, the \emph{scattering operator}. In
    general, $S=S(\lambda)$ depends on the eigenvalue (``energy'')
    parameter $\lambda$, namely, $S(\lambda)$ is (roughly) defined by
    looking how incoming and outgoing waves (of the form $x \mapsto
    \e^{\pm \im \sqrt \lambda x}$) propagate through a vertex.  In our
    case (i.e., if $L=0$ in $\laplacian{(\mc G,L)}$ described above),
    one can show that $S$ is independent of the \emph{energy}, namely,
    \begin{equation}
      \label{eq:scatt}
      S =
      \begin{pmatrix}
        \1 & 0 \\ 0 & -\1
      \end{pmatrix}
      = 2P - \1
    \end{equation}
    with respect to the decomposition $\Gmax=\mc G \oplus \mc
    G^\orth$, and where $P$ is the orthogonal projection of $\mc G$ in
    $\Gmax$.
  \item As in the discrete case, we can consider $\laplacian[0]{\mc G}
    := \laplacian {\mc G}$ as the Laplacian on $0$-forms, and
    $\laplacian[1]{\mc G} := \de_{\mc G} \de_{\mc G}^*$ as the
    Laplacian on $1$-forms, and again, by supersymmetry, we have the
    spectral relation
    \begin{equation*} \spec {\laplacian[1] {\mc G}} \setminus \{0\} =
      \spec {\laplacian[0]{\mc G}} \setminus \{0\}.
    \end{equation*}
    For more details and more general exterior derivatives
    corresponding to the case $L \ge 0$, we refer
    to~\cite[Sec.~5]{post:pre07a}.
  \end{enumerate}
\end{remark}

Using the definition $A \le B$ iff $\dom \qf a \supset \dom \qf b$ and
$\qf a(f) \le \qf b(f)$ for all $f \in \dom \qf b$ where $\qf a$, $\qf
b$ are the quadratic forms associated to the self-adjoint (unbounded)
non-negative operators $A$ and $B$ (i.e., $\qf a(f):=\normsqr {A^{1/2}
  f}=\iprod {Af} f$ for $f \in \dom \qf a := \dom A^{1/2}$ and $f \in
\dom A$, respectively), we have the following simple observation:
\begin{lemma}
  \label{lem:dn.brack}
  Assume that $\mc G_1 \le \mc G_2$ are two vertex spaces, then
  $\laplacian{\mc G_2} \le \laplacian {\mc G_1}$.
\end{lemma}
\begin{proof}
  The assertion follows directly from the inclusion $\Sobx {\mc G_1} X
  \subset \Sobx {\mc G_2} X$ and the fact that the quadratic forms are
  given by $\qf d_i(f) := \normsqr[\Lsqr X] {\de f}$ with $\dom \qf
  d_i=\Sobx{\mc G_i} X$.
\end{proof}
If $X$ is compact, i.e., the underlying graph is finite, we have:
\begin{proposition}
  \label{prp:res.comp}
  Assume that $X$ is compact, then the spectrum of $\laplacian{\mc G}$
  is purely discrete, i.e., there is an infinite sequence
  $\{\lambda_k\}_k$ of eigenvalues where
  $\lambda_k=\lambda_k(\laplacian{\mc G})=\lambda_k(\mc G)$ denotes
  the $k$-th eigenvalue (repeated according to its multiplicity) and
  $\lambda_k \to \infty$ as $k \to \infty$.
\end{proposition}
\begin{proof}
  We have to show that the resolvent of $\laplacian{\mc G}$ is a
  compact operator. This assertion follows easily from the estimate
  $\laplacian{\mc G} \ge \laplacian \Gmax=\bigoplus_e \laplacianN
  {I_e}$ where $\laplacianN {I_e}$ is the Neumann Laplacian on the
  interval $I_e$ having discrete spectrum
  $\lambda_k(\Gmax)=(k-1)^2\pi^2/\ell_e^2$ ($k=1,2,\dots$): The
  inequality implies the opposite inequality for the resolvents in
  $-1$; and therefore
  \begin{equation*}
    0 \le (\laplacian{\mc G}+1)^{-1} 
      \le (\laplacian \Gmax+1)^{-1}
        = \bigoplus_{e \in E} (\laplacianN {I_e}+1)^{-1}.
  \end{equation*}
  Since $E$ is finite, the latter operator is compact and therefore
  also the resolvent of $\laplacian{\mc G}$.
\end{proof}
Combining the last two results together with the variational
characterisation of the eigenvalues (the min-max principle), we have
the inequality
\begin{equation*}
  \lambda_k(\laplacian{\mc G_2}) \le \lambda_k(\laplacian{\mc G_1})
\end{equation*}
for all $k \in \N$ where $\mc G_1 \le \mc G_2$ are two vertex
spaces. Moreover,
\begin{equation*}
  \lambda_k^\Neu\bigl(\bigdcup_e I_e\bigr) = \lambda_k(\laplacian \Gmax) 
  \le \lambda_k(\laplacian {\mc G}) 
  \le \lambda_k(\laplacian{\mc G^{\min}}) 
   = \lambda_k^\Dir \bigl(\bigdcup_e I_e\bigr)
\end{equation*}
where $\lambda_k^\Dir \bigl(\bigdcup_e I_e\bigr)$ is the spectrum of
the (decoupled) Dirichlet operator $\laplacian{\mc
  G^{\min}}=\bigoplus_e \laplacianD{I_e}$. Note that $\lambda_k^\Neu
\bigl(\bigdcup_e I_e\bigr)=0$ for $k=1,\dots, \card E$, and
$\lambda_{k+\card E}^\Neu \bigl(\bigdcup_e I_e\bigr)=\lambda_k^\Dir
\bigl(\bigdcup_e I_e\bigr)$ where the latter sequence is a reordering
of the individual Dirichlet eigenvalues
$\lambda_k^\Dir(I_e)=k^2\pi^2/\ell_e^2$ repeated according to
multiplicity. In particular, for an equilateral metric graph (i.e,
$\ell_e=1$ for all edges $e$), then
\begin{equation*}
  (m-1)^2 \pi^2 \le \lambda_k(\laplacian {\mc G}) \le m^2 \pi^2,
  \qquad k=(m-1)\card E + 1,\dots, m \card E, \; m=1,2,\dots
\end{equation*}

For non-compact metric graphs, we can characterise the spectrum via
\emph{generalised} eigenfunctions, i.e., functions $\map f X \C$
satisfying the local vertex conditions $\ul f(v) \in \mc G_v$ and
$\orul f(v) \in \mc G^\orth_v$, but no integrability condition at
infinity: A measure $\rho$ on $\R$ is a \emph{spectral measure} for
$\laplacian{\mc G}$ iff for all measurable $I \subset \R$ we have
$\rho(I)=0$ iff the spectral projector satisfies $\1_I(\laplacian {\mc
  G})=0$.  In this case, we have the following result
(cf.~\cite[App.~B]{hislop-post:pre06}):
\begin{proposition}
  \label{prp:gen.ef}
  Assume the lower lengths bound~\eqref{eq:len.bd}. Let $\map \Phi X
  {(0,\infty)}$ be a bounded weight function, which is also in $\Lsqr
  X$. Then for almost every $\lambda \in \spec {\laplacian {\mc G}}$
  (with respect to a spectral measure), there is a generalised
  eigenfunction $f=f_\lambda$ associated to $\lambda$ such that
  \begin{equation*}
    \normsqr {\Phi f} = \int_X \abs {f(x)}^2 \Phi(x)^2 \dd x  < \infty.
  \end{equation*}
\end{proposition}
The function $\Phi$ can be constructed according to the graph. Denote
by $B_X(x_0,r)$ the metric ball of radius $r>0$ around the point $x_0
\in X$. For example, on a graph with sub-exponential volume growth,
i.e., for any $\eps>0$ there exists $C_\eps>0$ such that
\begin{equation*}
  \vol_1 B_X(x_0,r) := \int_X \1_{B_X(x_0,r)} \dd x
  \le C_\eps \e^{\eps r},
\end{equation*}
\sloppy the weight function $\Phi$ can be chosen in such a way that it
decays slower than exponentially, i.e., $\min \Phi(B_X(x_0,r)) \ge \wt
C_\eps \e^{-\eps r}$. In particular, we can choose $\Phi(x):=
\e^{-\eps d(x,x_0)} \le 1$, and, by Fubini, $\normsqr \Phi$ equals
\begin{equation*}
   \int_0^1 \vol_1\set{x \in X}{\Phi(x)^2 > t} \dd t
   = \int_0^1 \vol_1 B_X \Bigr(x_0, \frac {-\log t}{2\eps} \Bigr) \dd t
   \le C_\eps \int_0^1 t^{-1/2} \dd t < \infty.
\end{equation*}

\section{Relations between discrete and metric graphs}
\label{sec:rel.dg.mg}

In this section, we describe two cases, in which (parts of the)
spectrum of a metric graph can be described in terms of the discrete
graph. The first case deals with so-called \emph{equilateral} metric
graphs, i.e., graphs where all lengths are the same, say,
$\ell_e=1$. The second case treats the spectrum at the bottom, also in
the \emph{general} (non-equilateral) case.

\subsection{Equilateral metric graphs}
\label{sec:eql.mg}
An effective way of describing the relation between metric graph
Laplacians and the underlying (generalised) discrete one are so-called
\emph{boundary triples}. We do not give the general definition here.
instead, we refer to \cite{post:pre07c,bgp:08} and the references
therein. In brief, a boundary triple (originally developped for PDE
boundary value problems) describes an abstraction of Green's formula.

In order to describe the notions needed here, we define a
\emph{maximal} Laplacian in $\HS:= \Lsqr X$ with domain
\begin{equation*}
  \dom \laplacian[\max]{\mc G} := \Sobx[2] {\mc G} X
     = \bigset{f \in \Sobx[2]\max X} {\ul f \in \mc G},
\end{equation*}
i.e., we only fix the vertex values $\ul f$ to be in the vertex space
$\mc G$ with associated projection $P$. One can show similarly as in
the proof of \Prp{sa} that $\dom \laplacian[\max] {\mc G}$ is a closed
operator.

We define the boundary operators on the domain of the maximal
Laplacian as
\begin{subequations}
  \label{eq:bd.op}
  \begin{align}
    \map{\Gamma_0&} {\Sobx[2] {\mc G} X} {\mc G}, \qquad
    f \mapsto \ul f\\
    \map{\Gamma_1&} {\Sobx[2] {\mc G} X} {\mc G}, \qquad
    f \mapsto P \orul f'.
  \end{align}
\end{subequations}
Green's formula in this setting reads as
\begin{equation*}
  \iprod[\HS]{\laplacian[\max]{\mc G}f} g -
  \iprod[\HS] f {\laplacian[\max]{\mc G}g} =
  \iprod[\mc G]{\Gamma_0 f} {\Gamma_1 g} -
  \iprod[\mc G]{\Gamma_1 f} {\Gamma_0 g}
\end{equation*}
as one can easily see with the help of~\eqref{eq:part.int}.  As
self-adjoint reference operator, we denote by $\laplacian 0$ the
restriction of $\laplacian[\max]{\mc G}$ to $\ker \Gamma_0$. Note that
$\laplacian 0$ is precisely the metric graph Laplacian associated to
the minimal vertex space $\mc G^{\min}=0$, and therefore decoupled,
i.e.,
\begin{equation*}
  \laplacian 0 = \bigoplus_{e \in E} \laplacianD{I_e},
\end{equation*}
where $\laplacianD{I_e}$ denotes the Laplacian on $I_e$ with Dirichlet
boundary conditions and spectrum given by
$\spec{\laplacianD{I_e}}=\set {(\pi k/\ell_e)^2}{k =1,2,\dots}$ and
$\spec{\laplacian 0}$ is the union of all these spectra.

In the general theory of boundary triples, one can show that
$\Gamma_0$ restricted to ${\mc N^z=\ker(\laplacian[\max] {\mc G}-z)})$
is a topological isomorphism between $\mc N^z$ and $\mc G$ provided $z
\notin \spec {\laplacian 0}=:\Sigma$. We denote its inverse by
$\map{\beta(z)}{\mc G}{\mc N^z \subset \Lsqr X}$ (\emph{Krein's
  $\Gamma$-field}). In other words, $f=\beta(z)F$ is the solution of
the Dirichlet problem
\begin{equation*}
  (\lapl - z)f=0, \qquad
  \ul f = F.
\end{equation*}
Here, we can give an explicit formula for $\beta(z)$, namely
we have
\begin{equation*}
  f_e(x) = F_e(\bd_-e) s_{-,e,z}(x) + F_e(\bd_+e) s_{+,e,z}(x),
\end{equation*}
where\footnote{For $z=0$, we set $s_{-,e,0}(x):= 1-x/\ell_e$ and
  $s_{+,e,0}(x):=x/\ell_e$.}
\begin{equation}
  \label{eq:fund.sol}
  s_{-,e,z}(x) 
  = \frac {\sin (\sqrt z(\ell_e-x))}{\sin {\sqrt z \ell_e}}
  \und
  s_{+,e,z}(x) 
  = \frac {\sin (\sqrt z x)}{\sin {\sqrt z \ell_e}}.
\end{equation}
denote the fundamental solutions for $z \notin \spec{\lapl[0]}$.

Taking the derivative of $f=\beta(z)F$ on $\mc G$, i.e., defining
\begin{equation*}
  Q(z) F:= \Gamma_1 \beta(z)F,
\end{equation*}
we obtain a (bounded) operator $\map {Q(z)}{\mc G}{\mc G}$, called
\emph{Krein's Q-function} or \emph{Dirichlet-to-Neumann map}. Here, a
simple calculation shows that
\begin{equation*}
  (Q(z) F)_e(v) 
  = \frac {\sqrt z}{\sin (\sqrt z \ell_e)}
  \bigl[ \cos (\sqrt z \ell_e) F_e(v) - F_e(v_e) \bigr] .
\end{equation*}
if $z \notin \Sigma$.  In particular, if the metric graph is
equilateral (without loss of generality, $\ell_e=1$), we have
\begin{equation*}
  Q(z)
  = \frac {\sqrt z} {\sin \sqrt z}
  \bigl[ \dlapl[\mc G] - (1- \cos \sqrt z) \bigr].
\end{equation*}

The abstract theory of boundary triples gives here the following
result between the metric and discrete Laplacian. For a proof and more
general self-adjoint Laplacians as in \Remenum{sa}{other.sa} we refer
to~\cite[Sec.~5]{post:pre07c}. Certain special cases can be found for
example in~\cite{cattaneo:97,pankrashkin:06a,bgp:08}; and Pankrashkin
announced a more general result in~\cite{pankrashkin.talk:07}. For a
related result concerning a slightly different definition of a metric
graph Laplacian, see~\cite{baker-faber:06} and the references therein.
For spectral relations concerning averaging operators we refer
to~\cite{cartwright-woess:pre05}.
\begin{theorem}
  \label{thm:krein.mg}
  Assume the lower bound on the edge lengths~\eqref{eq:len.bd}.
  \begin{enumerate}
  \item
    \label{kernel}
    For $z \notin \spec{\lapl[0]}$ we have the explicit formula for
    the eigenspaces
    \begin{equation*}
      \ker (\laplacian{\mc G} - z) 
        = \beta(z) \ker Q(z).
      \end{equation*}
    \item
      \label{krein}
      For $z \notin \spec{\laplacian{\mc G}} \cup
      \spec{\lapl[0]}$ we have $0 \notin \spec{Q(z)}$ and Krein's
      resolvent formula
      \begin{equation*}
        (\laplacian{\mc G} - z)^{-1}
        = (\laplacian 0 - z)^{-1} -\beta(z) Q(z)^{-1} (\beta(\conj z))^*
      \end{equation*}
      holds.

    \item Assume that the graph is equilateral (say, $\ell_e=1$), then
      for $\lambda \in \C \setminus \R$ or $\lambda \in \R$ in the
      spectral gap $(\pi^2 k^2, \pi^2 (k+1)^2)$ ($k=1,2,\dots$) of
      $\laplacian 0$ or $\lambda < \pi^2$, we have
      \begin{equation*}
        (\laplacian{\mc G} - \lambda)^{-1}
        = (\laplacian 0 - \lambda)^{-1}
         -  \frac{\sin \sqrt \lambda}{\sqrt \lambda}
             \beta(\lambda) \bigl(\dlaplacian{\mc G}-(1-\cos \sqrt \lambda)
             \bigr)^{-1} (\beta(\conj \lambda))^*
      \end{equation*}
      and
      \begin{equation*}
        \lambda \in \spec[\bullet] 
                    {\laplacian{\mc G}} 
             \; \Leftrightarrow \;
               (1 - \cos \sqrt \lambda) \in 
                 \spec[\bullet] {\dlaplacian{\mc G}}
      \end{equation*}
      for all spectral types, namely, $\bullet \in \{\emptyset,
      \mathrm{pp}, \mathrm{disc}, \mathrm{ess} , \mathrm {ac}, \mathrm
      {sc}, \mathrm p\}$, the entire, pure point (set of all
      eigenvalues), discrete, essential, absolutely and singular
      continuous, and point spectrum ($\spec[p] A = \clo{\spec[pp]
        A}$). The multiplicity of an eigenspace is preserved.
  \end{enumerate}
\end{theorem}

\begin{remark}
  \label{rem:krein.qg}
  \begin{enumerate}
  \item The eigenspaces in \Thmenum{krein.mg}{kernel} for an
    equilateral graph can be constructed from the discrete data $F \in
    \ker(\dlaplacian{\mc G}-(1-\cos\sqrt z))$ by applying Krein's
    $\Gamma$-function, the ``solution operator'', namely, $f =
    \beta(z) F$ is the corresponding eigenfunction of the metric graph
    Laplacian.  The converse is also true: Given $f \in
    \ker(\laplacian{(\mc G, 0)}-z)$, then the corresponding
    eigenfunction $F \in \ker(\dlaplacian{\mc G}-(1-\cos\sqrt z))$ is
    just the restriction of $f$ to the vertices, namely $F=\ul f$.
  \item The resolvent formula in \Thmenum{krein.mg}{krein} is very
    explicit, since
    \begin{equation*}
      (\lapl[0] - z)^{-1} = \bigoplus_{e \in E} (\laplacianD{I_e}-z)^{-1}
    \end{equation*}
    is decoupled and explicit formulas for the resolvent on the
    interval are known.  In particular, the analysis of the
    (equilateral) metric graph resolvent is reduced to the analysis of
    the discrete Laplacian resolvent (see
    also~\cite{kostrykin-schrader:06,kps:pre07}).

    Krein's resolvent formula~\eqref{krein} is very useful when
    analysing further properties of the quantum graph $(X,
    \laplacian{\mc G})$ via the resolvent.
  \item We excluded the Dirichlet spectrum $\spec {\laplacian 0} =
    \Sigma$. These values may occur in the spectrum of $\laplacian
    {\mc G}$ or not. For example, if $\mc G$ is the standard vertex
    space $\mc G^\stand$ and if $X$ contains a loop with an even
    number of edges each having the same length, we can define on each
    edge a Dirichlet solution on the edge (with opposite sign on
    successive edges). This function is continuous in the vertices,
    and satisfies also the Kirchhoff condition in each vertex.
    Therefore, on a metric graph, compactly supported eigenfunctions
    may exist.
  \end{enumerate}
\end{remark}

\subsection{Relation at the bottom of the spectrum}
\label{sec:ind.discr}
%
Let us analyse the spectrum at the bottom in more detail.  As in
\Sec{vx.sp} we define the Hilbert chain associated to the exterior
derivative $\de_{\mc G}$ as
\begin{equation*}
   \mc C_{X,\mc G} \colon 0 
   \longrightarrow \Sobx {\mc G} X \stackrel {\de_\mc G} 
   \longrightarrow \Lsqr X \longrightarrow 0
\end{equation*}
and call elements of the first space \emph{$0$-forms}, and of the
second space \emph{$1$-forms}.  The associated cohomology spaces (with
coefficients in $\C$) are defined as
    \begin{align*}
      H^0(X,\mc G) &:= \ker \de_{\mc G} \cong \ker \de_{\mc G}/\ran 0, \\
      H^1(X,\mc G) &:= \ker \de^*_{\mc G} = \ran \de_{\mc G}^\orth
      \cong \ker 0/\ran \de_{\mc G}
    \end{align*}

The \emph{index} or \emph{Euler characteristic} of the cohomology
associated to the metric graph $X$ with vertex space $\mc G$ is then
defined as
\begin{equation*}
  \ind (X,\mc G) := \dim \ker \de_{\mc G} - \dim \ker \de_{\mc G}^*,
\end{equation*}
i.e., the \emph{Fredholm index} of $\de_{\mc G}$, provided at least
one of the dimensions is finite.

We have the following result (for more general cases
cf.~\cite{post:pre07a}, and for a different approach
see~\cite{fkw:07}):
\begin{theorem}
  \label{thm:index}
  Assume that $G$ is a weighted discrete graph with lower lengths
  bound~\eqref{eq:len.bd}, and denote by $X$ the associated metric
  graph, and by $\mc G$ a vertex space associated to $G$. Then there
  is an isomorphism $\Phi^*=\Phi_0^* \oplus \Phi_1^*$ with
  \begin{equation*}
    \map {\Phi_p^*} {H^p(X,\mc G)} {H^p(G,\mc G)}.
  \end{equation*}
  More precisely, $\Phi^*$ is induced by a Hilbert chain morphism
  $\Phi$ , i.e.,
  \begin{equation*}
    \begin{diagram}
      \mc C_{X,\mc G} \colon 
           0 & \rTo & \Sobx {\mc G} X & \rTo^{\de_{\mc G}} &
      \Lsqr X & \rTo & 0\\
      & & \dTo{\Phi_0} & &
      \dTo{\Phi_1}               &      &     \\
      \mc C_{G,\mc G} \colon 
           0 & \rTo & \mc G & \rTo^{\dde_{\mc G}} & \lsqr E & \rTo & 0
    \end{diagram}
  \end{equation*}
  is commutative, where
  \begin{equation*}
    \Phi_0 f := \ul f = \Gamma_0 f, \qquad
    \Phi_1 g := \Bigl\{\int_{I_e} g_e(x) \dd x\Bigr\}_e.
  \end{equation*}
  In particular, if $G$ is finite (and therefore $X$ compact), then
  \begin{equation*}
    \ind(G,\mc G)=\ind(X,\mc G).
  \end{equation*}
\end{theorem}
For general results on Hilbert chains and their morphisms we refer
to~\cite[Ch.~1]{lueck:02} or~\cite{bruening-lesch:92}.
\begin{proof}
  The operators $\Phi_p$ are bounded. Moreover, that $\Phi$
  is a chain morphism follows from
  \begin{equation*}
    (\Phi_1 \de_{\mc G} f)_e = \int_{I_e} f'_e(x) \dd x
    = f_e(\ell_e)-f_e(0) = (\dde_{\mc G} \ul f)_e 
    = (\dde_{\mc G} \Phi_0 f)_e.
  \end{equation*}
  Furthermore, there is a Hilbert chain morphism $\Psi$, i.e., 
  \begin{equation*}
    \begin{diagram}
      \mc C_{X,\mc G} \colon 
      0 & \rTo & \Sobx {\mc G} X & \rTo^{\de_{\mc G}} &
      \Lsqr X & \rTo & 0\\
      & & \uTo{\Psi_0} & &
      \uTo{\Psi_1}               &      &     \\
     \mc C_{G,\mc G} \colon 
      0 & \rTo & \mc G & \rTo^{\dde_{\mc G}} & \lsqr E & \rTo & 0
    \end{diagram}
  \end{equation*}
  given by
  \begin{equation*}
    \Psi_0 F := \beta(0) F 
    = \{ F_e(\bd_-e)s_{-,e,0} + F_e(\bd_+e) s_{+,e,0} \}_e, \qquad
    \Psi_1 \eta := \{ \eta_e \1_{I_e}/\ell_e \}_e
  \end{equation*}
  (see~\Eq{fund.sol}), i.e., we let $\Phi_0 F$ be the affine
  (harmonic) function on $I_e$ with boundary values fixed; and $\Phi_1
  \eta$ be an (edgewise) constant function. Again, the chain morphism
  property $\Psi_1 \dde_{\mc G}= \de_{\mc G} \Psi_0$ can easily be
  seen. Furthermore, $\Phi \Psi$ is the identity on the second
  (discrete) Hilbert chain $\mc C_{G,\mc G}$.  It follows now from
  abstract arguments (see e.g.~\cite[Lem.~2.9]{bruening-lesch:92})
  that the corresponding induced maps $\Phi_p^*$ are isomorphisms on
  the cohomology spaces.
\end{proof}

\begin{remark}
  The sub-complex $\Psi(\mc C_{G,\mc G})$ of $\mc C_{X,\mc G}$ consists
  of the subspace of edge-wise affine functions ($0$-forms) and of
  edge-wise constant functions ($1$-forms). In this way, we can
  naturally embed the discrete setting into the metric graph one.  In
  particular, assume that $0<\ell_0 \le \ell_e \le \ell_+ < \infty$,
  then
  \begin{multline*}
    \normsqr{\Psi_0 F}
     = \sum_e \frac 1 {\ell_e} \int_0^{\ell_e} 
          \abs{  F_e(\bd_-e) (\ell_e - x)
               + F_e(\bd_+e) x }^2 \dd x\\
     = \sum_e \frac 1 {3\ell_e} \int_0^{\ell_e} 
          \abs{  F_e(\bd_-e)^2 + F_e(\bd_-e)F_e(\bd_+e)
               + F_e(\bd_+e)^2},
  \end{multline*}
  so that
  \begin{equation*}
    \frac 1{2\ell_+} \normsqr[\mc G] F \le
      \normsqr{\Psi_0 F} \le  \frac 5 {6\ell_0} \, \normsqr[\mc G] F,
  \end{equation*}
  i.e., redefining the norm on $\mc G$ by $\norm[\mc G,1] F
  :=\norm{\Psi_0 F}$ gives an equivalent norm turning $\Psi_0$ into an
  isometry.  Moreover, $\norm{\Psi_1 \eta} = \norm[\lsqr E] \eta$.
  For more details on this point of view (as well as ``mixed'' types
  of discrete and metric graphs), we refer
  to~\cite{friedman-tillich:04} and references therein.
\end{remark}
\section{Relations between metric graphs and manifolds}
\label{sec:mfd}

\begin{figure}[h]
  \begin{center}
\begin{picture}(0,0)%
 \includegraphics{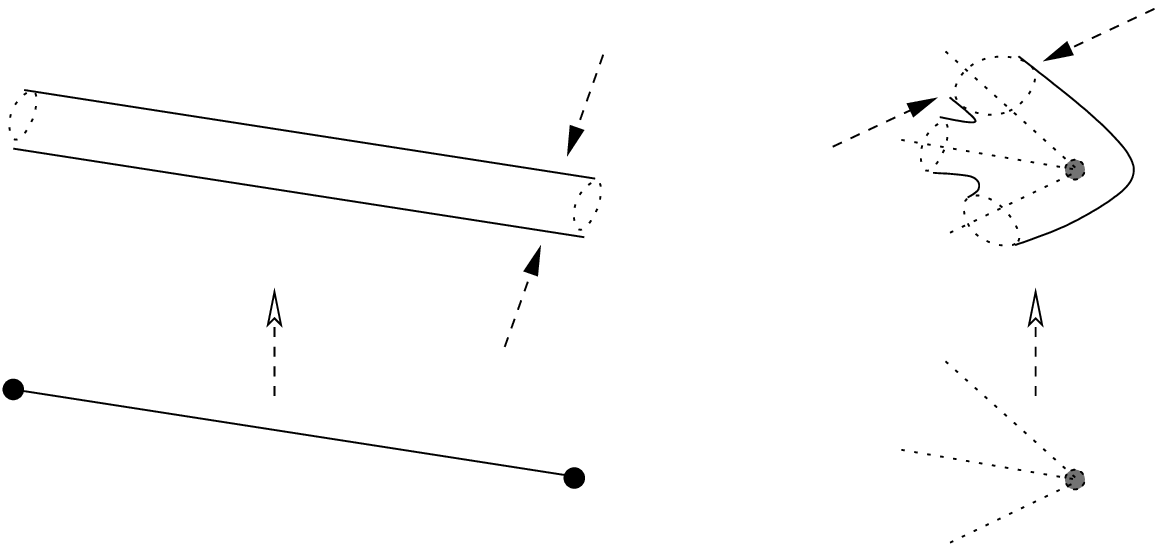}%
\end{picture}%
\setlength{\unitlength}{4144sp}%
\begin{picture}(5282,2500)(374,-1970)
  \put(1591,-586){$U_\edeps$}%
  \put(1561,-1674){$I_e$}%
  \put(5641,-331){$U_\vxeps$}%
  \put(5416,-1771){$v$}%
  \put(4579,374){$\eps$}%
  \put(3169,-436){$\eps$}%
\end{picture}
    \caption{The associated edge and vertex neighbourhoods with
      $F_\eps=\Sphere^1_\eps$, i.e., $U_\edeps$ and $U_\vxeps$ are
      $2$-dimensional manifolds with boundary.}
    \label{fig:edge.vertex}
  \end{center}
\end{figure}

Let us briefly describe the relation of a metric graph
$X_0=(V,E,\bd,\ell)$ with manifolds.  For more details, we refer to
the review article~\cite{exner-post:pre07b} and the references
therein.  Let $X_\eps$ be a $d$-dimensional connected manifold with
metric $g_\eps$.  If $X_\eps$ has boundary, we denote it by $\bd
X_\eps$; let us stress that our discussion covers different kind of
models, like the $\eps$-neighbourhood of an metric graph
\emph{embedded} in $\R^\nu$, as well as sleeve-type manifolds (like
the surface of a pipeline network) having no boundary. We assume that
$X_\eps$ can be decomposed into open sets $U_\edeps$ and $U_\vxeps$,
i.e.,\footnote{
  \label{fn:dcup} The expression $A=\bigdcup_i A_i$ means that the
  $A_i$'s are open (in $A$), mutually disjoint and the interior of
  $\bigcup_i \clo A_i$ equals $A$.}
\begin{equation}
  \label{eq:decomp}
  X_\eps = \bigdcup_{e \in E} U_\edeps \,\dcup\,
  \bigdcup_{v \in V} U_\vxeps.
\end{equation}
Denote the metric on $X_\eps$ by $g_\eps$.  To simplify the discussion
here, we assume that $U_\edeps$ and $U_\vxeps$ are isometric to
\begin{subequations}
  \label{eq:met}
  \begin{align}
    U_\edeps & \cong (I_e \times F, g_\edeps) &
    g_\edeps &= \de x_e^2 + \eps^2 h\\
    U_\vxeps & \cong (U_v, g_\vxeps) & g_\vxeps &= \eps^2 g_v
  \end{align}
\end{subequations}
where $(F,h)$ is a compact $m$-dimensional manifold with $m:=(d-1)$,
and $(U_v,g_v)$ is an $\eps$-independent $d$-dimensional manifold
(cf.~\Fig{edge.vertex}). Strictly speaking, for a metric graph $X_0$
embedded in $\R^2$, the edge neighbourhoods of the associated
$\eps$-neighbourhood $X_\eps := \set {x \in \R^2}{\dist(x,X_0) <
\eps}$ must be shorter in the longitudinal direction in order to have
space for the vertex neighbourhoods. Nevertheless, this fact causes
only an error of order $\eps$ for the associated metrics, which does
not matter in our convergence analysis below.

Note that $\bd U_\vxeps \setminus \bd X_\eps$ has $(\deg v)$-many
components isometric to $(F,\eps^2 h)$ denoted by $(\bd_e U_v, \eps^2
h)$ for $e \in E_v$.  The cross section manifold $F$ has a boundary or
does not have one, depending on the analogous property of $X_\eps$.

On the other hand, given a metric graph $X_0$ and vertex neighbourhood
manifolds $U_v$ as below, we can abstractly construct a graph-like
manifold $X_1$ from these building blocks according to the rules of
the graph with a family of metrics $g_\eps$
satisfying~\eqref{eq:decomp} and~\eqref{eq:met} with
$X_\eps=(X_1,g_\eps)$.

For simplicity, we suppose that $\vol_m F=1$.  Then we have
\begin{equation}
  \label{eq:met.prod.asym}
  \dd U_\edeps = \eps^m \dd F \dd x_e
\end{equation}
for the Riemannian densities.  We consider the Hilbert space $\HS_\eps
= \Lsqr {X_\eps}$ and the Laplacian $\laplacian {X_\eps} := \de^* \de
\ge 0$ (with Neumann boundary conditions if $\bd X_\eps \ne
\emptyset$) where $\de$ denotes the exterior derivative.
In addition, we assume that the following uniformity conditions
are valid,
\begin{align}
  \label{eq:vol.ew}
  c_{\vol} := \sup_{v \in V}\: {\vol_d U_v} < \infty, \qquad \qquad
  \lambda_2 := \inf_{v \in V} \EWN 2 {U_v} > 0,
\end{align}
where $\EWN 2 {U_v}$ denotes the second (i.e., first non-zero) Neumann
eigenvalue of $(U_v,g_v)$.  Roughly speaking, the
requirements~\eqref{eq:vol.ew} mean that the region $U_v$ remains
small w.r.t. the vertex index. Obviously, these assumptions are
trivially satisfied once the vertex set $V$ is finite.

In order to compare operators in $\Lsqr {X_0}$ and $\Lsqr
{X_\eps}=\bigoplus_e \Lsqr {I_e} \otimes \Lsqr {F_\eps} \oplus
\bigoplus_v \Lsqr {U_\vxeps}$, we use the identification operator
\begin{equation}
  \label{eq:j}
  J f := \{ f_e \otimes \1_\eps \}_e \oplus \{ 0_v\}_v,
\end{equation}
where $\1_\eps=\eps^{-m/2}\1$ is the lowest normalised eigenfunction
on $F_\eps=(F,\eps^2h)$ and in turn $0_v$ is the zero function on
$U_v$.  This identification operator is \emph{quasi-unitary}, i.e.,
$J^*J=\id_{\HS_0}$, $\norm J = 1$, and one can show that
\begin{equation*}
  \norm{(JJ^*-\id_{\HS_\eps})(\laplacian{X_\eps}+1)^{-1/2})} 
   = \Err(\eps^{1/2}),
\end{equation*}
where $\Err(\eps^{1/2})$ depends only on the lower lengths bound
$\ell_0$ in~\eqref{eq:len.bd}, on $c_{\vol}$ and $\lambda_2$:
Basically, we have
\begin{equation*}
    (\id_{\HS_\eps}-JJ^*) u 
   = \sum_{e \in E}\int_{I_e} 
         \normsqr[F_\eps]{ u(x,\cdot)- 
              \iprod[F_\e\eps]{u(x,\cdot)}{\1_\eps}} \dd x
 + \sum_{v \in V} \normsqr[U_\vxeps] u
\end{equation*}
and both contributions can be estimated in terms of $\Err(\eps)
(\normsqr{\de u} + \normsqr u)$ as stated above.

The main statement of this section is the following
(cf.~\cite{post:06}):
\begin{theorem}
  \label{thm:res}
  Assume the uniformity conditions~\eqref{eq:len.bd}, $d_0 := \sup_v
  \deg v < \infty$ and~\eqref{eq:vol.ew}. Then the Laplacians
  $\laplacian {X_\eps}$ and $\laplacian {X_0}$ are
  $\Err(\eps^{1/2})$-close with respect to the quasi-unitary map $J$
  defined in~\eqref{eq:j}, i.e.
  \begin{equation*}
    \norm{(\laplacian {X_\eps}+1)^{-1} J
          - J (\laplacian {X_0} + 1)^{-1}}
    \le \Err(\eps^{1/2}).
  \end{equation*}
   In addition, we have
  \begin{equation}
    \label{eq:sandwich}
    \norm{(\laplacian {X_\eps}+1)^{-1}
          - J (\laplacian {X_0} + 1)^{-1} J^*}
    \le \Err(\eps^{1/2}),
  \end{equation}
  where the error term depends only on $\ell_0$, $d_0$, $c_{\vol}$ and
  $\lambda_2$.
\end{theorem}
Once the above resolvent estimates are established, one can develop a
functional calculus for pairs of operators
$(\HS_\eps,\laplacian{X_\eps})$ and $(\HS_0,\laplacian{X_0})$ together
with a quasi-unitary identification operator $J$, and establish the
above operator estimates also for more general functions $\phi$ than
$\phi(\lambda)=(\lambda+1)^{-1}$, namely for spectral projectors
($\phi=\1_I$, $I$ interval), or for the heat operator
($\phi_t(\lambda)=\e^{-t \lambda}$, $t > 0$). Moreover, we can show
that the spectra are close to each other:
\begin{theorem}
  \label{thm:res.spec}
  Under the assumptions of the previous theorem, we have
  \begin{equation*}
    \norm{ \1_I(\laplacian {X_\eps}) J
          - J \1_I(\laplacian {X_0})}
    \le \Err(\eps^{1/2})
      \quad \text{and} \quad
    \norm{ \1_I(\laplacian {X_\eps})
          - J \1_I(\laplacian {X_0})J^*}
    \le \Err(\eps^{1/2})
  \end{equation*}
  for the spectral projections provided $I$ is a compact interval such
  that $\bd I \cap \spec {\laplacian {X_0}} = \emptyset$. In
  particular, if $I$ contains a single eigenvalue $\lambda(0)$ of
  $\laplacian {X_0}$ with multiplicity one corresponding to an
  eigenfunction $u(0)$, then
  there is an eigenvalue $\lambda(\eps)$ and an eigenfunction
  $u(\eps)$ of $\laplacian{X_\eps}$ such that
  \begin{equation*}
    \norm{J u(0) - u(\eps)} = \Err(\eps^{1/2}).
  \end{equation*}
  In addition, the spectra converge uniformly on $[0,\Lambda]$, i.e.
  \begin{equation*}
    \spec {\laplacian {X_\eps}} \cap [0,\Lambda] \to
    \spec {\laplacian {X_0}} \cap [0,\Lambda]
  \end{equation*}
  in the sense of the Hausdorff distance on compact subsets of
  $[0,\Lambda]$. The same result is true if we consider only the
  essential or the discrete spectral components.
\end{theorem}
In particular, the above theorem applies to the case when $X_0$ (and
therefore $X_\eps$) is compact, and we obtain
\begin{equation}
  \label{eq:ew}
  \lambda_k(\laplacian {X_\eps}) - \lambda_k (\laplacian {X_0})
   = \Err(\eps^{1/2}).
\end{equation}
This estimate can also be proved directly by applying the min-max
theorem, and estimating the errors of the corresponding Rayleigh
quotients. For more results (like a similar convergence of resonances)
we refer again to~\cite{exner-post:pre07b} and the references therein.
Recently, Grieser showed in~\cite{grieser:pre07} an asymptotic
expansion of the eigenvalues and the eigenfunctions also for other
boundary conditions on $\bd X_\eps$, for example Dirichlet.

\section{Estimates on the first non-zero eigenvalue}
\label{sec:cheeger}
Here, we comment on inequalities on the first non-zero eigenvalue of a
graph, namely a lower bound in terms of an isoperimetric constant. For
details, see e.g.~\cite{chung:97,higuchi-shirai:04,nicaise:87}.

Let $X$ be a compact metric graph and $Y \subset X$ be a non-empty
open subset. We denote by $\card {\bd Y}$ the number of points in the
boundary (``volume'' of dimension $0$), and by $\vol_1 Y:=\int_X \1_Y
\dd x$ the total length of $Y$ (``volume'' of dimension $1$).
\emph{Cheeger's (isoperimetric) constant} for the metric graph $X$ is
defined as
\begin{equation}
  \label{eq:cheeger}
  h(X):= \inf_Y \frac {\card {\bd Y}}
              {\min(\vol_1 Y,\vol_1 \compl Y)}
\end{equation}
where $\compl Y := X \setminus Y$, and the infimum runs over all open,
subset $Y \subset X$ such  that $Y \ne \emptyset$ and $Y \ne X$.

For simplicity, we assume that $X$ is connected and that each vertex
space is standard, i.e., $\mc G_v=\mc G^\stand = \C(1,\dots,1)$. The
corresponding (standard or Kirchhoff) Laplacian (denoted by
$\laplacian[\stand] X$) has discrete spectrum. In particular, the
first eigenvalue fulfills $\lambda_1(\laplacian[\stand] X)=0$, while
the second is positive $\lambda_2(X) := \lambda_2(\laplacian[\stand]
X) > 0$. If not already obvious, this follows from \Thm{index}, and
the fact that the dimension of the $0$-th cohomology group counts the
number of components.

Cheeger's theorem in this context is the following:
\begin{theorem}
  \label{thm:cheeger}
  Assume that $X$ is a connected, compact metric graph with standard
  vertex space $\mc G^\stand$ and denote by $\lambda_2(X) > 0$ the
  first non-zero eigenvalue of the standard (Kirchhoff) metric graph
  Laplacian. Then we have
  \begin{equation*}
    \lambda_2(X) \ge \frac {h(X)^2}4.
  \end{equation*}
\end{theorem}
\begin{proof}
  The proof follows closely the line of arguments as in the manifold
  case (see also \cite{nicaise:87}). The basic ingredient is the
  \emph{co-area formula}
  \begin{equation*}
    \int_X \abs{\phi'(x)} \dd x 
    = \int_0^\infty \card{\set{x \in X}{\phi(x)=t}} \dd t
  \end{equation*}
  for any non-negative, edgewise $\Contsymb^1$-function $\phi$.

  Denote by $f$ the corresponding eigenfunction associated to
  $\lambda_2(X)$. Without loss of generality, we may assume that $f$
  is real-valued. Set $X_+:=\set{x \in X}{f(x) > 0}$.  Moreover, we
  may assume that $\vol_1 X_+ \le \vol_1 \compl{X_+}$ (if this is not
  true, replace $f$ by $-f$). Finally, $X_+ \ne \emptyset$ and $X_+
  \ne X$ since $f$ changes sign as second eigenfunction (only the
  first eigenfunction is constant).

  Let $g := \1_{X_+} f$, then $g$ is non-negative and $g \ne 0$.
  Moreover, since $g$ is continuous, we can perform partial
  integration without additional boundary terms in $\bd X_+$: In
  particular, if $v \in \bd X_+$ is a vertex then $g_e(v)=0$ for all
  adjacent edges $e \in E_v$. In particular, we have
  \begin{equation}
    \label{eq:proof.ch}
    \lambda_2(X) = \frac {\iprod g {-g''}} {\normsqr g}
    = \frac {\normsqr{g'}} {\normsqr g}
    \ge \frac 1 4 \Bigg( \frac {\int_X \abs{(g^2)'(x)} \dd x}
                               {\normsqr g} \Bigg)^2
  \end{equation}
  where we used Cauchy-Schwarz for the latter inequality. Setting
  $X(t):= \set{x \in X}{g(x)^2 > t}$, the co-area formula and Fubini
  yield
  \begin{equation*}
    \frac {\int_X \abs{(g^2)'(x)} \dd x} {\normsqr g}
    = \frac {\int_0^\infty \card{\set{x \in X}{g(x)^2=t}} \dd t}
            {\int_0^\infty \vol_1 X(t) \dd t}
    \ge \frac {\int_0^{t_0} \card{\bd X(t)} \dd t}
            {\int_0^{t_0} \vol_1 X(t) \dd t}
  \end{equation*}
  since $\{g^2=t\} \supset \bd X(t)$. Here, $t_0 := \max g(X)^2 > 0$
  because $X_+ \ne \emptyset$.  Moreover, $X(t)$ is open ($g$ is
  continuous), $\vol_1 X(t) \le \vol_1 X_+ \le \vol_1 \compl X_+ \le
  \vol_1 \compl{X(t)}$, $X(t) \ne \emptyset$ for $t \in [0,t_0)$ and
  $X(t) \ne X$ for all $t \ge 0$ since $X_+ \ne X$. The definition of
  Cheeger's constant finally yields the lower bound $h(X)$ for the
  last fraction.
\end{proof}
\begin{remark}
  One might ask whether similar results hold for more general boundary
  spaces $\mc G$ (i.e., the metric graph Laplacian $\laplacian{\mc
    G}$). There are several problems in the general case:
  \begin{itemize}
  \item If the projection $P_v$ associated to $\mc G_v$ has complex
    entries, the eigenfunction may no longer be chosen to  be
    real-valued.
  \item If the function $f \in \dom \laplacian {\mc G}$ is not
    continuous at a vertex, (e.g., negative on one edge and positive
    on another edge meeting in the same vertex), the boundary terms of
    $g$ appearing from partial integration in~\eqref{eq:proof.ch} may
    not vanish at this vertex.
  \item The eigenfunction associated to the first non-zero eigenvalue
    may not change its sign (e.g., if it is a Dirichlet function on a
    single edge). In this case, one needs a modified Cheeger constant
    (with $\vol_1 Y$ in the denominator, and $Y \subset X$ open, not
    intersecting the ``Dirichlet'' vertices.
  \end{itemize}
\end{remark}
Cheeger's theorem for a (standard) finite discrete graph $G=(V,E,\bd)$
can be proven in a similar way. For simplicity, we assume that all
weights are the same, say $\ell_e=1$, and that the graph has no
self-loops. We define Cheeger's constant for the discrete graph $G$ as
\begin{equation}
  \label{eq:cheeger.dg}
  h(G):= \inf_W \frac {\card {E(W,\compl W)}}
              {\min(\vol_0 W,\vol_0 \compl W)},
\end{equation}
where the infimum runs over all subsets $W \subset V$ such that $W \ne
\emptyset$ and $W \ne V$. Furthermore, $E(W,\compl W)$ is the set of
all edges having one vertex in $W$ and the other one in $\compl W$.
The volume of $W$ is defined as $\vol_0 W := \sum_{v \in W} \deg w$.
Note that $\vol_0 W = \normsqr[\lsqr V]{\1_W}$
(see~\eqref{eq:norm.std}).  For a proof of the next theorem, see
e.g.~\cite[Thm.~2.2]{chung:97}.
\begin{theorem}
  \label{thm:cheeger.dg}
  Assume that $G$ is a connected, finite discrete graph with standard
  vertex space $\mc G^\stand$ and denote by $\lambda_2(G) > 0$ the
  first non-zero eigenvalue of the standard discrete graph Laplacian
  as defined in~\eqref{eq:lap.std}. Then we have
  \begin{equation*}
    \lambda_2(G) \ge \frac {h(G)^2}2.
  \end{equation*}
\end{theorem}
Again, it would be interesting to carry over the above result for more
general discrete Laplacians, namely for $\dlaplacian{\mc G}$ and a
general vertex space $\mc G$ associated to $G$.

Let us finally mention an upper bound on the second eigenvalue in
terms of the distance of subsets (see~\cite{friedman-tillich:04}
or~\cite{cgy:96} for the general scheme and a similar result for
discrete graphs):
\begin{theorem}
  \label{thm:upper}
  Let $X$ be a connected, compact metric graph and denote by
  $\lambda_2(X)$ the second (first non-zero) eigenvalue of the
  standard metric graph Laplacian on $X$. Then
  \begin{equation*}
    \lambda_2(X) \le \frac 4 {d(A,B)^2} \Bigg(\log \frac {\vol_1
      X}{(\vol_1 A \vol_1 B)^{1/2}} \Bigg)^2
  \end{equation*}
  for any two disjoint measurable subsets $A,B$ of $X$, where $d(A,B)$
  denotes the distance between the sets $A$ and $B$ in the metric
  graph $X$.
\end{theorem}
One can prove similar results also for higher eigenvalues.
Note that \Thms{cheeger}{upper} also hold in the manifold case (with
the appropriate measures), and that they are consistent with the
eigenvalue approximation result of~\eqref{eq:ew}.

\section{Trace formulas}
\label{sec:trace}
In this last section we present some results concerning the trace of
the heat operator. Trace formulas for metric graph Laplacians appeared
first in an article of Roth~\cite{roth:84} (see
also~\cite{kurasov:07}), where he used standard (Kirchhoff) boundary
conditions; more general self-adjoint vertex conditions
(energy-independent, see \Remenum{sa}{scatt}) are treated
in~\cite{kostrykin-schrader:06,kps:pre07}.

We first need some (technical) notation; inevitable in order to
properly write down the trace formula. For simplicity, we assume that
the graph has no self-loops.

\begin{definition}
  A \emph{combinatorial path} in the discrete graph $G$ is a sequence
  $c=(e_0,v_0,e_1,v_1,\dots,e_n,v_n,e_{n+1})$ where $v_i \in \bd e_i
  \cap \bd e_{i+1}$ for $i=0, \dots, n$. We call $\card c := n+1$ the
  \emph{combinatorial length} of the path $c$, and $e_-(c):=e_0$
  resp.\ $e_+(c):=e_{n+1}$ the \emph{initial} resp.\ \emph{terminal
    edge} of $c$.
  Similarly, we denote by $\bd_- c:= v_0$ and $\bd_+c := v_n$ the
  \emph{initial} resp.\ \emph{terminal vertex} of $c$, i.e., the first
  resp.\ last vertex in the sequence $c$.  A \emph{closed path} is a
  path where $e_-(c)=e_+(c)$. A closed path is \emph{properly closed}
  if $c$ is closed and $\bd_- c \ne \bd_+ c$. Denote by $C_n$ the set
  of all properly closed paths of combinatorial length $n$, and by $C$
  the set of all properly closed paths.
\end{definition}
If the graph does not have double edges, a properly closed
combinatorial path can equivalently be described by the sequence
$c=(v_0, \dots, v_n)$ of vertices passed by. In particular, $\card
{C_0} = \card V$, $\card {C_1}=0$ (no self-loops) and $\card
{C_2}=2\card E$. Moreover, $C_3=\emptyset$ is equivalent that $G$ is
bipartite. A graph $G$ is called \emph{bipartite}, if $V=V_+ \dcup
V_-$ with $E=E(V_+,V_-)$, see \Eq{cheeger.dg}.

\begin{definition}
  Two properly closed paths $c$, $c'$ are called \emph{equivalent} if
  they can be obtained from each other by successive application of
  the cyclic transformation
  \begin{equation*}
    (e_0,v_0,e_1,v_1,\dots,e_n,v_n,e_0) \to
    (e_1,v_1,\dots,e_n,v_n,e_0, v_0, e_1.)
  \end{equation*}
  The corresponding equivalence class is called \emph{cycle} and is
  denoted by $\wt c$. The set of all cycles is denoted by $\wt C$.
  Given $p \in \N$ and a cycle $\wt c$, denote by $p \wt c$ the cycle
  obtained from $\wt c$ by repeating it $p$-times. A cycle $\wt c$ is
  called \emph{prime}, if $\wt c = p \wt c'$ for any other cycle $\wt
  c'$ implies $p=1$. The set of all prime cycles is denoted by $\wt
  C_\prim$.
\end{definition}

\begin{definition}
  Let $\map \gamma {[0,1]} X$ be a \emph{metric path} in the metric
  graph $X$, i.e., a continuous function which is of class
  $\Contsymb^1$ on each edge and $\gamma'(t) \ne 0$ for all $t \in
  [0,1]$ such that $\gamma(t) \in X^1=X \setminus V$, i.e., inside an
  edge. In particular, a path in $X$ cannot turn its direction inside
  an edge. We denote the set of all paths from $x$ to $y$ by
  $\Gamma(x,y)$.
\end{definition}
Associated to a metric path $\gamma \in \Gamma(x,y)$ there is a unique
combinatorial path
$c_\gamma=(e_0,v_0,e_1,v_1,\dots,e_{n-1},v_n,e_{n+1})$ determined by
the sequence of edges and vertices passed along $\gamma(t)$ for $0 < t
< 1$, (it is not excluded that $\gamma(0)$ or $\gamma(1)$ is a vertex;
this vertex is not encoded in the sequence $c$). In particular, if
$x=\gamma(0), y=\gamma(1) \notin V$, then $x$ is on the initial
edge $e_-(c)$ and $y$ on the terminal edge $e_+(c)$.

On the other hand, a combinatorial path $c$ and two points $x,y$ being
on the initial resp.\ terminal edge, i.e., $x \in \clo e_-(c)$, $y \in
\clo e_+(c)$, but different from the initial resp.\ terminal vertex,
i.e., $x \ne \bd_-(c)$ and $y \ne \bd_+(c)$, uniquely determine a
metric path $\gamma=\gamma_c \in \Gamma(x,y)$ (up to a change of
velocity).  Denote the set of such combinatorial paths from $x$ to $y$
by $C(x,y)$.

\begin{definition}
  The \emph{length} of the metric path $\gamma \in \Gamma(x,y)$ is
  defined as $\ell(\gamma):= \int_0^1 \abs{\gamma'(s)} \dd s$. In
  particular, if $c=c_\gamma=(e_0,v_0, \dots, e_n,v_n,e_{n+1})$ is the
  combinatorial path associated to $\gamma$, then
  \begin{equation*}
    d_c(x,y):= \ell(\gamma) = 
    \abs{x-\bd_- c_\gamma}
    + \sum_{i=1}^n \ell_{e_i}
    + \abs{y -\bd_- c_\gamma},
  \end{equation*}
  where $\abs{x-y}:=\abs{x_e-y_e}$ denotes the distance of $x,y$ being
  inside the same edge $e$ (or its closure), and $x_e$, $y_e \in \clo
  I_e$ are the corresponding coordinates ($x=\Phi_e x_e$, cf.\
  \Remenum{met.space}{coord}). Note that there might be a shorter path
  between $x$ and $y$ \emph{outside} the edge $e$.  For a properly
  closed path $c$ we define the \emph{metric length} of $c$ as
  $\ell(c)=\ell(\gamma_c)$ and similarly, $\ell(\wt c):= \ell(c)$ for
  a cycle. Note that this definition is well-defined.
\end{definition}
Finally, we need to define the \emph{scattering amplitudes}
associated to a combinatorial path $c=(e_0,v_0, \dots,
e_n,v_n,e_{n+1})$ and a vertex space $\mc G$. Denote by $P=\oplus P_v$
its orthogonal projection in $\Gmax$ onto $\mc G$. Denote by $S:=
2P-\1$ the corresponding scattering matrix defined in \Eq{scatt}. In
particular, $S$ is local, i.e., $S=\oplus_v S_v$ and we define
\begin{equation*}
  S_{\mc G}(c) := \prod_{i=0}^n S_{e_i,e_{i+1}}(v_i),
\end{equation*}
where $S_{e,e'}(v)=2P_{e,e'}(v)-\delta_{e,e'}$ for $e,e' \in E_v$.
For a cycle, we set $S(\wt c):= S(c)$, and this definition is
obviously well-defined, since multiplication of complex numbers is
commutative.  

For example, the standard vertex space $\mc G^\stand$ has projection
$P=(\deg v)^{-1} \mathbb E$ (all entries are the same), so that
\begin{align*}
  S^\stand_{e,e'}(v)&=\frac 2 {\deg v}, \quad e \ne e',&
  S^\stand_{e,e}(v)&= \frac 2 {\deg v} - 1.
\end{align*}
If in addition, the graph is regular, i.e, $\deg v=r$ for all $v \in
V$, then one can simplify the scattering amplitude of a combinatorial
path $c$ to
\begin{equation*}
  S^\stand(c) = \Bigl(\frac 2 r \Bigl)^a \Bigl(\frac 2 r -1 \Bigl)^b 
\end{equation*}
where $b$ is the number of reflections in $c$ ($e_i=e_{i+1}$) and $a$
the number of transmissions $e_i\ne e_{i+1}$) in $c$.

We can now formulate the trace formula for a compact metric graph with
Laplacian $\laplacian{\mc G}$ (cf.~\cite[Thm.~1]{roth:84},
\cite[Thm.~4.1]{kps:pre07}):
\begin{theorem}
  \label{thm:trace}
  Assume that $X$ is a compact metric graph (without self-loops), $\mc
  G$ a vertex space and $\laplacian{\mc G}$ the associated
  self-adjoint Laplacian (cf.\ \Prp{sa}). Then we have
  \begin{equation*}
    \tr \e^{-t \laplacian {\mc G}} =
    \frac {\vol_1 X}{2(\pi t)^{1/2}}
    + \frac 12 \bigl( \dim \mc G - \card E \bigr) + 
    \frac 1{2(\pi t)^{1/2}} 
    \sum_{\wt c \in \wt C_\prim} \sum_{p \in \N}
    S_{\mc G}(\wt c)^p \ell(\wt c) 
        \exp\Bigl(-\frac {p^2 \ell(\wt c)^2}{4t}\Bigr)
  \end{equation*}
  for $t > 0$, where $\vol_1 X = \sum_e \ell_e$ is the total length of
  the metric graph $X$.
\end{theorem}
\begin{remark}
  \indent
  \begin{enumerate}
  \item The first term in the RHS is the term expected from the Weyl
    asymptotics. The second term is precisely $1/2$ of the index $\ind
    (X,\mc G)$ of the metric (or discrete) graph $X$ with vertex space
    $\mc G$, i.e., the Fredholm index of $\de_{\mc G}$. In \Thm{index}
    we showed that the index is the same as the discrete index $\ind
    (G,\mc G)$ (the Fredholm index of $\dde_{\mc G)}$).
    In~\cite{kps:pre07}, the authors calculated the second term as
    $(\tr S)/4$, but since $S=2P-\1$, we have $\tr S=2\dim\mc G-\dim
    \Gmax=2(\dim \mc G - \card E)$. The last term in the trace formula
    comes from an combinatorial expansion.
  \item The sum over prime cycles of the metric graph $X$ is an
    analogue of the sum over primitive periodic geodesics on a
    manifold in the celebrated Selberg trace formula, as well as an
    analogue of a similar formula for (standard) discrete graphs, see
    \Thm{trace.dg}.
  \item Trace formulas can be used to solve the inverse problem: For
    example, Gutkin, Smilansky and Kurasov,
    Nowaczyk~\cite{gutkin-smilansky:01,kurasov:07,kurasov-nowaczyk:05}
    showed that if $X$ does not have self-loops and double edges, and
    if all its lengths are rationally independent, then the metric
    structure of the graph is uniquely determined. Further extensions
    are given e.g.  in~\cite{kps:pre07}. Counterexamples
    in~\cite{roth:84,gutkin-smilansky:01,bss:06} show that the
    rational independence is really needed, i.e., there are
    isospectral, non-homeomorphic graphs.
  \end{enumerate}
\end{remark}
The proof of \Thm{trace} uses the expansion of the heat kernel, namely
one can show that
\begin{equation*}
  p_t(x,y) = \frac 1 {2 (\pi t)^{1/2}} 
    \Bigl( \delta_{x,y} \exp\Bigl(-\frac {\abs{x-y}^2}{4t}\Bigr) 
       + \sum_{c \in C(x,y)} S(c) 
             \exp\Bigl(-\frac {d_c(x,y)^2}{4t}\Bigr)\Bigr),
\end{equation*}
where $\delta_{x,y}=1$ if $x,y$ are inside the closure of the same
edge (and not both on opposite sides of $\bd e$) and $0$ otherwise.
The trace of $\e^{-t \laplacian{\mc G}}$ can now be calculated as the
integral over $p_t(x,x)$. The first term in the heat kernel expansion
gives the volume term, the second splits into \emph{properly} closed
paths leading to the third term (the sum over prime cycles), and the
index term in the trace formula is the contribution of non-properly
closed paths.  More precisely, a non-properly closed path runs through
its initial and terminal edge (which are the same by definition of a
closed path) in opposite directions. For more details, we refer to
\cite{roth:84} or \cite{kps:pre07}.

Let us finish with some trace formulas for discrete graphs. Assume for
simplicity, that $G$ is a simple discrete graph, i.e., $G$ has no
self-loops and double edges, and that all lengths are the same
($\ell_e=1$). For simplicity, we write $v \sim w$ if $v,w$ are
connected by an edge. Let $\mc G$ be an associated vertex space. Since
$\laplacian {\mc G}=\1-M_{\mc G}$ and $M_{\mc G}$ (see~\Eq{markov})
are bounded operators on $\mc G$ , we have
\begin{equation*}
  \tr \e^{-t \laplacian {\mc G}} 
  = \e^{-t} \tr \e^{t M_{\mc G}}
  = \e^{-t} \sum_{n=0}^\infty \frac {t^n}{n!} \tr M_{\mc G}^n.
\end{equation*}
Furthermore, using~\eqref{eq:markov} $n$-times, we obtain
\begin{equation*}
  M_{\mc G}^n=\bigoplus_{v_0} \sum_{v_1 \sim v} \dots
               \sum_{v_n \sim v_{n-1}} 
   A_{\mc G}(v_0,v_1) A_{\mc G}(v_1,v_2) \cdot \ldots \cdot
         A_{\mc G}(v_{n-1},v_n),
\end{equation*}
and
\begin{equation*}
  \tr M_{\mc G}^n= \sum_{v_0} \sum_{v_1 \sim v_0} \dots
               \sum_{v_{n-1} \sim v_{n-2}} 
   \tr A_{\mc G}(v_0,v_1) A_{\mc G}(v_1,v_2) \cdot \ldots \cdot
         A_{\mc G}(v_{n-1},v_0).
\end{equation*}
Note that the sum is precisely over all combinatorial, (properly)
closed paths $c=(v_0, \dots, v_{n-1}) \in C_n$. Denoting by
\begin{equation*}
  W_{\mc G}(c) := 
   \tr A_{\mc G}(v_0,v_1) A_{\mc G}(v_1,v_2) \cdot \ldots \cdot
        A_{\mc G}(v_{n-1},v_0).
\end{equation*}
the \emph{weight} associated to the path $c$ and the vertex space $\mc
G$, we obtain the following general trace formula. In particular, we
can write the trace as a (discrete) ``path integral'':
\begin{theorem}
  \label{thm:trace.dg}
  Assume that $G$ is a discrete, finite graph with weights $\ell_e=1$
  having no self-loops or double edges. Then
  \begin{equation}
    \label{eq:trace.vx}
    \tr \e^{-t \dlaplacian {\mc G}} 
    = \e^{-t} \sum_{n=0}^\infty  \sum_{c \in C_n} 
    \frac {t^n}{n!} W_{\mc G}(c) 
    = \e^{-t} \sum_{c \in C} \frac {t^{\card c}}{\card c!} 
    W_{\mc G}(c).
  \end{equation}
\end{theorem}
Let us interprete the weight in the standard case $\mc G=\mc
G^\stand$. Here, $A_{\mc G^\stand}(v,w)$ can be interpreted as
operator from $\C(\deg w)$ to $\C(\deg v)$ (the degree indicating the
corresponding $\lsymb_2$-weight) with $A_{\mc G^\stand}(v,w)=1$ if
$v,w$ are connected and $0$ otherwise. Viewed as multiplication in
$\C$ (without weight), $A_{\mc G^\stand}(v,w)$ is unitarily equivalent
to the multiplication with $(\deg v \deg w)^{-1/2}$ if $v \sim w$
resp.\ $0$ otherwise. In particular, if $c=(v_0, \dots, v_{n-1})$ is
of length $n$, then the weight is
\begin{equation*}
  W^\stand(c) = \frac 1 {\deg v_0} \cdot \frac 1 {\deg v_1} \cdot \ldots \cdot
             \frac 1 {\deg v_{n-1}}.
\end{equation*}
If, in addition, $G$ is a regular graph, i.e., $\deg v=r$ for all $v
\in V$, then $W^\stand(c)=r^{-n}$. Then the trace
formula~\eqref{eq:trace.vx} reads as
\begin{equation*}
  \tr \e^{-t \laplacian {\mc G^\stand}} 
  = \e^{-t} \sum_{n=0}^\infty
         \frac {t^n}{r^n n!} \card{C_n}
         = \e^{-t} \Bigr( \card V + \frac {\card E} {2 r^2} \, t^2 +
           \frac {\card{C_3}}{6 r^3} \, t^3 + \dots \Bigr),
\end{equation*}
since $\card {C_0}=\card V$, $\card {C_1}=0$ (no self-loops) and
$\card {C_2}=2\card E$. In particular, one can determine the
coefficients $\card {C_n}$ form the trace formula expansion.

The weight $W^\stand(c)$ for the standard vertex space is a sort of
probability of a particle chosing the path $c$ (with equal probability
to go in any adjacent edge at each vertex). It would be interesting to
give a similar meaning to the ``weights'' $W_{\mc G}(c)$ for general
vertex spaces.

\def\cprime{$'$} \def\cprime{$'$} \def\cprime{$'$} \def\cprime{$'$}
  \def\cprime{$'$} \def\cprime{$'$} \def\cprime{$'$} \def\cprime{$'$}
\providecommand{\bysame}{\leavevmode\hbox to3em{\hrulefill}\thinspace}
\renewcommand{\MR}[1]{}

\providecommand{\MRhref}[2]{%
  \href{http://www.ams.org/mathscinet-getitem?mr=#1}{#2}
}
\providecommand{\href}[2]{#2}

\end{document}